\documentclass[aop,seceqn,citesort,MSNbibl,dvips]{arximspdf}

\doi{10.1214/10-AOP598}
\volume{39}
\issue{6}
\pubyear{2011}
\firstpage{2385}
\lastpage{2423}

\makeatletter

\newtheorem{Proposition}{Proposition}
\newtheorem{Theorem}[Proposition]{Theorem}

\newcommand{\p}{+}
\newcommand{\m}{-}

\newcommand{\EE}{\mathsf E}
\newcommand{\PP}{\mathsf P}
\newcommand{\cF}{\mathcal{F}}
\newcommand{\cB}{\mathcal{B}}
\newcommand{\R}{\mathbb R}
\newcommand{\LL}{\mathbb L}
\newcommand{\eps}{\varepsilon}

\makeatother

\begin{document}
\begin{frontmatter}

\title{Predicting the ultimate supremum of a stable L\'evy process with no negative jumps}
\runtitle{Predicting the ultimate supremum}

\begin{aug}
\author[A]{\fnms{Violetta} \snm{Bernyk}\ead[label=e1]{violetta.bernyk@ubs.com}},
\author[B]{\fnms{Robert C.} \snm{Dalang}\ead[label=e2]{robert.dalang@epfl.ch}} and
\author[C]{\fnms{Goran} \snm{Peskir}\corref{}\ead[label=e3]{goran@maths.man.ac.uk}}
\runauthor{V. Bernyk, R. C. Dalang and G. Peskir}
\affiliation{UBS AG, Ecole Polytechnique F\'ed\'erale and The
University of Manchester}
\address[A]{V. Bernyk\\
UBS AG \\
Europastrasse 1 \\
8152 Opfikon \\
Switzerland \\
\printead{e1}} 
\address[B]{R. C. Dalang\\
Institut de Math\'ematiques \\
Ecole Polytechnique F\'ed\'erale \\
Station 8 \\
1015 Lausanne \\
Switzerland \\
\printead{e2}}
\address[C]{G. Peskir \\
School of Mathematics \\
The University of Manchester \\
Oxford Road \\
Manchester M13 9PL \\
United Kingdom \\
\printead{e3}}
\end{aug}

\received{\smonth{1} \syear{2010}}

%
\begin{abstract}
Given a stable L\'evy process $X=(X_t)_{0 \le t \le T}$ of index
$\alpha\in(1,2)$ with no negative jumps, and letting $S_t =
\sup_{ 0 \le s \le t} X_s$ denote its running supremum for $t \in
[0,T]$, we consider the optimal prediction problem
\[
V = \inf_{0 \le\tau\le T} \EE(S_T \m X_\tau)^p,
\]
where the infimum is taken over all stopping times $\tau$ of $X$,
and the error parameter $p \in(1,\alpha)$ is given and fixed.
Reducing the optimal prediction problem to a fractional
free-boundary problem of Riemann--Liouville type, and finding an
explicit solution to the latter, we show that there exists $\alpha_*
\in(1,2)$ (equal to $1.57$ approximately) and a strictly increasing
function $p_* \dvtx (\alpha_*,2) \rightarrow(1,2)$ satisfying
$p_*(\alpha_*+) = 1$, $p_*(2-) = 2$ and $p_*(\alpha) < \alpha$ for
$\alpha\in(\alpha_*,2)$ such that for every $\alpha\in
(\alpha_*,2)$ and $p \in(1,p_*(\alpha))$ the following stopping
time is optimal
\[
\tau_* = \inf\{ t \in[0,T] \dvtx S_t \m X_t \ge z_* (T \m
t)^{1/\alpha} \},
\]
where $z_* \in(0,\infty)$ is the unique root to a transcendental
equation (with parameters $\alpha$ and $p$). Moreover, if either
$\alpha\in(1,\alpha_*)$ or $p \in(p_*(\alpha),\alpha)$ then it is
not optimal to stop at $t \in[0,T)$ when $S_t \m X_t$ is
sufficiently large. The existence of the breakdown points $\alpha_*$
and~$p_*(\alpha)$ stands in sharp contrast with the Brownian motion
case (formally corresponding to $\alpha=2$), and the phenomenon
itself may be attributed to the interplay between the jump structure
(admitting a~transition from lighter to heavier tails) and the
individual preferences (represented by the error parameter $p$).
\end{abstract}

%

\pdfkeywords{60G40, 60J75, 45J05, 60G25, 47G20, 26A33, Optimal
prediction, optimal stopping, ultimate supremum, stable Levy process with no negative
jumps, spectrally positive, fractional free-boundary problem, Riemann--Liouville fractional
derivative, Caputo fractional derivative, stochastic process reflected at its
supremum, infinitesimal generator, weakly singular Volterra integral
equation, polar kernel, smooth fit, curved boundary}

\end{frontmatter}

\section{Introduction}

Stopping a stochastic process $X=(X_t)_{0 \le t \le T}$ as close as
possible to its ultimate supremum $S_T = \sup_{ 0 \le s \le T} X_s$
is an objective of both practical and theoretical interest. Speaking
in general terms, the optimal prediction problem can be formulated
as follows
%
%
\begin{equation}
\label{1.1}
V = \inf_{0 \le\tau\le T} d(X_\tau,S_T),
\end{equation}
where the infimum is taken over all stopping times $\tau$ of $X$,
and $d$ is a~distance/error function [e.g., $d(X_\tau,S_T)=
\EE(S_T \m X_\tau)^p$ where $p>0$ is a parameter quantifying the
error]. Variants of these problems have been studied in the past
mostly in discrete time (see, e.g., \cite{Ka,GM,Bo,GS}), and the case
of continuous time has been studied in the
recent papers \cite{GPS} and \cite{Pe} when~$X$ is a standard
Brownian motion. This study was extended in~\cite{DP-1} to the case
of Brownian motion with drift. It was observed there that the
existence of a~nonzero drift leads to optimal stopping boundaries
having a complex structure which in some cases appears to be
counter-intuitive. For other optimal prediction problems studied to
date, we refer to \cite{Sh-1,Ur,DP-2,Sh-2,DPS} (see also~\cite{PS},
Chapter VIII). In these problems, it
is assumed that the underlying process has continuous sample paths.

The purpose of the present paper is to initiate a study of the
optimal prediction problems for processes with jumps in continuous
time, and to examine the extent to which the jump structure
influences the resulting optimal stopping boundaries. To stay close
to the more familiar case of Brownian motion, we study the case when
$X$ is a stable L\'evy process of index $\alpha\in(1,2)$, and to
focus on one particular aspect of the jump structure we consider the
case when $X$ jumps upward only (i.e., when $X$ has no negative
jumps). It turns out that already these hypotheses lead to a
complicated optimal prediction problem, which apart from initial
similarities with the case of Brownian motion (through the scaling
property and deterministic time-change arguments) requires novel
arguments to be developed in order to find a solution. These
complications are primarily attributed to the underlying jump
structure which leads to the relatively unexplored avenue of
integro-differential equations (fractional calculus) instead of more
familiar differential equations. Yet another difficulty (that the
law of $S_T$ was not available in the literature prior to the
present study) is now overcome by the accompanying paper \cite{BDP},
and the knowledge of this law plays a key role in our treatment of
the optimal prediction problem below.

Our main findings (Theorem \ref{theo11}) can be summarized as follows. Given a
stable L\'evy process $X=(X_t)_{0 \le t \le T}$ of index $\alpha\in
(1,2)$ with no negative jumps, and letting $S_t = \sup_{ 0 \le s
\le t} X_s$ denote its\vspace*{1pt} running supremum for $t \in[0,T]$, we
consider the optimal prediction problem
%
%
\begin{equation}
\label{1.2}
V = \inf_{0 \le\tau\le T} \EE(S_T \m X_\tau)^p,
\end{equation}
where the infimum is taken over all stopping times $\tau$ of $X$,
and the error parameter $p \in(1,\alpha)$ is given and fixed (we
will see in Section \ref{sec2} below why the restriction to this interval is
natural). Reducing the optimal prediction problem to a fractional
free-boundary problem of Riemann--Liouville type, and finding an
explicit solution to the latter, we show that there exists $\alpha_*
\in(1,2)$ (equal to $1.57$ approximately) and a strictly increasing
function $p_* \dvtx (\alpha_*,2) \rightarrow(1,2)$
satisfying $p_*(\alpha_*+)=1, p_*(2-) = 2$ and $p_*(\alpha) <
\alpha$ for $\alpha\in(\alpha_*,2)$ such that for every $\alpha
\in(\alpha_*,2)$ and $p \in(1,p_*(\alpha))$ the following stopping
time is optimal
%
%
\begin{equation}
\label{1.3}
\tau_* = \inf\{ t \in[0,T] \dvtx S_t \m X_t \ge z_* (T \m
t)^{1/\alpha} \},
\end{equation}
where $z_* \in(0,\infty)$ is the unique root to a transcendental
equation (with parameters $\alpha$ and $p$). This extends the
analogous results for a standard Brownian motion $X$ derived in
\cite{GPS} and \cite{Pe} when $p=2$ and $p \in(1,2)$, respectively.
Moreover, if either $\alpha\in(1,\alpha_*)$ or $p \in
(p_*(\alpha),\alpha)$ then it is not optimal to stop at $t \in
[0,T)$ when $S_t \m X_t$ is sufficiently large. The existence of the
breakdown points $\alpha_*$ and $p_*(\alpha)$ stands in sharp
contrast with the Brownian motion case (formally corresponding to
$\alpha=2$), and the phenomenon itself may be attributed to the
interplay between the jump structure (admitting a transition from
lighter to heavier tails) and the individual preferences
(represented by the error parameter $p$). In particular, recalling
that the index $\alpha$ quantifies the heaviness of the upward tails
of the process $X$, we see that the result may be broadly
interpreted as follows: \textit{the heavier the upward tails the
larger the optimal stopping time}. While this conclusion is close to
naive intuition, and the interpretation itself may also be extended
to account for the individual preferences, the fact that the
solution method can detect the breakdown points exactly appears to
be of considerable practical and theoretical interest. Other
interesting features of the problem include the remarkable
probabilistic representation of the solution to the
It\^o/Riemann--Liouville/Caputo free-boundary problem that is novel
in the case of Brownian motion as well.

\section{The optimal prediction problem}\label{sec2}

1. Let $X=(X_t)_{t \ge0}$ be a stable L\'evy process of index
$\alpha\in(1,2)$ whose characteristic function is given by
%
%
\begin{equation}
\label{2.1}
\EE e^{i \lambda X_t} = \exp\biggl(t \int_0^\infty(e^{i \lambda
x} - 1 - i \lambda x) \frac{c}{x^{1+\alpha}} \,dx\biggr) = e^{c
\Gamma(-\alpha) (-i \lambda)^\alpha t}
\end{equation}
for $\lambda\in\R$ and $t \ge0$ with $c>0$. Let $S=(S_t)_{t \ge
0}$ denote the supremum process of~$X$, that is,
%
%
\begin{equation}
\label{2.2}
S_t = \sup_{0 \le s \le t} X_s
\end{equation}
for $t \ge0$. Consider the optimal prediction problem
%
%
\begin{equation}
\label{2.3}
V = \inf_{0 \le\tau\le T} \EE(S_T \m X_\tau)^p,
\end{equation}
where the infimum is taken over all stopping times $\tau$ of $X$
[i.e., stopping times with respect to the natural filtration
$\cF_t^X = \sigma(X_s \dvtx 0 \le s \le t)$ generated by $X$ for $t \ge
0$]. It is assumed in (\ref{2.3}) that the error parameter $p \in
(1,\alpha)$ and the terminal time $T>0$ are given and fixed (we will
see below that there is no restriction in assuming that $T=1$).

2. The following properties of $X$ are readily deduced from
(\ref{2.1}) using standard means (see, e.g., \cite{Be} and
\cite{Ky}): the law of $(X_{\sigma t})_{t \ge0}$ is the same as the
law of $(\sigma^{1/\alpha} X_t)_{t \ge0}$ for each $\sigma>0$ given
and fixed (scaling property); $X$ is a martingale with $\EE X_t = 0$
for all $t \ge0$; $X$ jumps upward (only) and creeps downward
[in the sense that $\PP(X_{\rho_x} =x)=1$ for $x<0$ where $\rho_x
= \inf\{ t \ge0 \dvtx X_t < x \}$ is the first entry time of $X$
into $(-\infty,x)$]; $X$ has sample paths of unbounded variation;
$X$ oscillates from $-\infty$ to $+\infty$ (in the sense that
$\liminf_{ t \rightarrow\infty} X_t = -\infty$ and $\limsup_{ t
\rightarrow\infty} X_t = +\infty$ both a.s.); the starting point
$0$ of $X$ is regular [for both $(-\infty,0)$ and $(0,+\infty)$].
Note also that the L\'evy measure $\nu$ of $X$ equals
%
%
\begin{equation}
\label{2.4}
\nu(dx) = \frac{c}{x^{1+\alpha}} \,dx
\end{equation}
on the Borel $\sigma$-algebra of $(0,\infty)$. Setting, for example,
$c =
1/(2\Gamma(-\alpha))$ we see from (\ref{2.1}) that $X =$ $X(\alpha)$
converges in law to a standard Brownian motion $B$ as $\alpha
\uparrow2$. We moreover see from (\ref{2.4}) that when $\alpha$ is
closer to $2$ then the (upward) jumps of $X$ have lighter tails, and
when $\alpha$ is closer to $1$ then the (upward) jumps of $X$ have
heavier tails. Thus, in many ways, the process $X$ resembles a
standard Brownian motion $B$, however, the existence of (upward)
jumps of $X$ represents a notable exception. Note also that $X_t$ is
not equal in law to $-X_t$ for fixed $t>0$ unlike in the case of
$B$.

3. The error parameter $p$ in the problem (\ref{2.3}) is assumed to
belong to $(1,\alpha)$ for two reasons. First, it is well known
(see, e.g., \cite{Sa}, page 159) that for a L\'evy process $X=(X_t)_{t
\ge0}$ and a number $p>0$ given and fixed, the following three
facts are equivalent: (i)~$\EE X_t^p < \infty$ for some/all $t>0$;
(ii)~$\EE\sup_{ 0 \le s \le t} X_s^p < \infty$ for some/all
$t>0$; (iii)~$\int_1^\infty x^p \nu(dx) < \infty$. In the case
of our process $X$ when $\nu$ is given by (\ref{2.4}) above, it is
easily seen that (iii)~holds [and thus both expected values in (i)
and (ii) are finite] if and only if $p < \alpha$. In particular, the
latter condition then also implies that the value $V$ in (\ref{2.3})
is finite. Second, if $p=1$ then the optimal prediction problem
(\ref{2.3}) is trivial since $\EE X_\tau= 0$ for every (bounded)
stopping time $\tau$ of $X$ due to the martingale property of $X$.
Hence, $p \in(1,\alpha)$ represents a natural assumption on the
error parameter.

4. Note that there is no loss of generality if we assume that $T=1$
in the problem (\ref{2.3}). Indeed, if we set $V=V(T)$ to indicate
dependence on $T>0$ in (\ref{2.3}), then by the scaling property of
$X$ we see that $V(T) = T^{p/\alpha} V(1)$ and there is a simple
one-to-one correspondence between the stopping times $\tau$ in the
problem $V(T)$ and the stopping times $\sigma$ in the problem $V(1)$
(obtained by setting $\sigma=\tau/T$). For this reason, we will
often assume in the sequel that the horizon $T$ in (\ref{2.3})
equals $1$.

5. \textit{Projecting future onto present}. One of the key initial
difficulties in the optimal prediction problem (\ref{2.3}) is that
the expression after the expectation sign contains the random
variable $S_T$ and as such depends on the (ultimate) future of the
process $X$ that is unknown at the present (stopping) time $\tau\in
[0,T)$. In our first step therefore (similarly to \cite{GPS} and
\cite{Pe}), we will project the future states of $X$ onto the
present/past states of $X$ by conditioning with respect to
$\cF_\tau^X$ and exploiting stationary/independent increments of
$X$. As already mentioned above, we may and do assume that $T=1$ in
the sequel.\vadjust{\goodbreak}

To this end, note that we have
%
%
\begin{eqnarray}
\label{2.5}
&&\EE\bigl((S_1 \m X_t)^p \vert\cF_t^X\bigr)\nonumber\\
&&\qquad= \EE\Bigl( \Bigl(
\sup_{0 \le s \le t} (X_s \m X_t) \vee\sup_{t \le s \le1}
(X_s \m X_t) \Bigr)^{ p} \big\vert\cF_t^X \Bigr) \\
&&\qquad= \bigl(\EE(y \vee S_{1-t})^p \bigr)\big\vert_{y = S_t-X_t}\nonumber
\end{eqnarray}
since $\sup_{ t \le s \le1} (X_s \m X_t) \stackrel{\mathrm{law}}{=}
S_{1-t}$ is independent from $\cF_t^X$ and $S_t \m X_t$ is
$\cF_t^X$-measurable. Moreover, we can write
%
%
\begin{eqnarray}
\label{2.6}
\EE(y \vee S_{1-t})^p &=& \int_0^\infty\PP\bigl( (y \vee S_{1-t})^p
> z \bigr) \,dz \nonumber\\
&=& y^p + \int_{y^p}^\infty\PP( S_{1-t}^p > z
) \,dz \nonumber\\
&=& y^p + \int_{y^p}^\infty\PP\bigl( (1 \m
t)^{p/\alpha} S_1^p > z \bigr) \,dz \nonumber\\[-8pt]\\[-8pt]
&=& (1 \m t)^{p/\alpha}
\biggl[ \biggl( \frac{y}{(1 \m t)^{1/\alpha}} \biggr)^{ p}\nonumber\\
&&\hspace*{50.53pt}{} + \int_{(
{y}/{(1-t)^{1/\alpha}})^p}^\infty\PP( S_1^p > w )
\,dw \biggr] \nonumber\\
&=&\!: F(t,y)\nonumber
\end{eqnarray}
upon using that $S_{1-t} \stackrel{\mathrm{law}}{=} (1 \m
t)^{1/\alpha} S_1$ by the scaling property of $X$ and substituting
$w = z/(1 \m t)^{p/\alpha}$. Combining (\ref{2.5}) and (\ref{2.6}),
we get
%
%
\begin{equation}
\label{2.7}
\EE\bigl((S_1 \m X_t)^p \vert\cF_t^X\bigr) = F(t,S_t \m X_t)
\end{equation}
for all $t \ge0$. Using the fact that each stopping time $\tau$ of
$X$ is the limit of a~decreasing sequence of discrete stopping times
$\tau_n$ of $X$ as $n \rightarrow\infty$, it is easily verified
using Hunt's lemma (see, e.g., \cite{Wi}, page 236) that (\ref{2.7})
extends as follows
%
%
\begin{equation}
\label{2.8}
\EE\bigl((S_1 \m X_\tau)^p \vert\cF_\tau^X\bigr) = F(\tau,S_\tau
\m X_\tau)
\end{equation}
for all stopping times $\tau$ of $X$ with values in $[0,1]$. Setting
%
%
\begin{equation}
\label{2.9}
Y_t = S_t \m X_t
\end{equation}
for $t \ge0$ it is well known (see, e.g., \cite{Be}) that
$Y\,{=}\,(Y_t)_{t\ge0}$ is a time-homogeneous (strong) Markov process
with respect to $(\cF_t^X)_{t \ge0}$ (obtained by reflecting $X$ at
its supremum $S$). Taking $\EE$ on both sides in (\ref{2.8}) and
using the notation (\ref{2.9}), we see that the optimal prediction
problem (\ref{2.3}) reduces to the optimal stopping problem
%
%
\begin{equation}
\label{2.10}
V = \inf_{0 \le\tau\le1} \EE F(\tau,Y_\tau),
\end{equation}
where the infimum is taken over all stopping times $\tau$ of $X$.
This optimal stopping problem is two-dimensional (see, e.g.,
\cite{PS}, Section 6) since the underlying (strong) Markov process is
the time--space process $((t,Y_t))_{0 \le t \le1}$ and the
horizon~$1$ is finite. We will now show (similarly to \cite{GPS}) that this
problem can further be reduced to a one-dimensional infinite-horizon
optimal stopping problem for a (killed) Markov process $Z=(Z_s)_{s
\ge0}$. It should be noted that the time-change arguments used in
\cite{GPS} when $X$ is a standard Brownian motion are not directly
applicable in the present context (due to the absence of L\'evy's
characterization theorem).

5. \textit{Deterministic time change}. Motivated by the form of the
function $F$ in~(\ref{2.6}), we now introduce the deterministic time
change
%
%
\begin{equation}
\label{2.11}
t(s) = 1 \m e^{-\alpha s},
\end{equation}
where $t(s) \in[0,1)$ is the ``old'' time and $s \in[0,\infty)$ is a
``new'' time. Note that $\tau= t(\sigma)$ is a stopping time with
respect to $(\cF_t^X)_{t \ge0}$ if and only if $\sigma=
t^{(-1)}(\tau)$ is a stopping time with respect to
$(\cF_{t(s)}^X)_{s \ge0}$. Letting $F_{S_1}$ denote the
distribution function of $S_1$ and setting
%
%
\begin{equation}
\label{2.12}
G(z) = \EE(z \vee S_1)^p= z^p + \int_{z^p}^\infty\bigl( 1 -
F_{S_1}(w^{1/p}) \bigr) \,dw
\end{equation}
for $z \ge0$, we see from (\ref{2.6}) and (\ref{2.12}) that
%
%
\begin{equation}
\label{2.13}
F(t,S_t \m X_t) = e^{-ps} G(Z_s)
\end{equation}
for all $t=t(s) \in[0,1)$ and all $s \in[0,\infty)$ satisfying
(\ref{2.11}), where $Z=(Z_s)_{s \ge0}$ is a new stochastic process
defined by
%
%
\begin{equation}
\label{2.14}
Z_s = e^s \bigl(S_{t(s)} \m X_{t(s)} \bigr)
\end{equation}
for $s \ge0$. It turns out that $Z$ is a time-homogeneous (strong)
Markov process. Moreover, the following proposition reveals that one
can enable $Z$ to start at arbitrary points and still preserve the
(strong) Markov property. This fact will play a prominent role in
the main proof below.
\begin{Proposition}\label{prop1}
The stochastic process $Z=(Z_s)_{s \ge
0}$ defined in (\ref{2.14}) is a time-homogenous (strong) Markov
process with respect to the filtration $(\cF_{t(s)}^X)_{s \ge0}$.
Moreover, if we set
%
%
\begin{equation}
\label{2.15}
Z_s^z = e^s \bigl(z \vee S_{t(s)} \m X_{t(s)} \bigr)
\end{equation}
for $s \ge0$ and $z \in\R_+$, then $\PP_{ z} :=
\mathrm{Law}((Z_s^z)_{s \ge0} \vert\PP)$ defines a family of
probability measures on the canonical space of c\`adl\`ag functions
$(D_+,\cB(D_+))$ under which the coordinate process $C=(C_s)_{s \ge
0}$ is (strong) Markov with $\PP_{ z}(C_0 = z)=1$ for $z \in
\R_+$.
\end{Proposition}
\begin{pf}
We have
%
%
\begin{eqnarray}
\label{2.16}
Z_{s+h}^z &=& e^{s+h} \bigl(z \vee S_{t(s+h)} \m X_{t(s+h)} \bigr)
\nonumber\\
&=& e^{s+h} \Bigl(\Bigl[\bigl(z \vee S_{t(s)} \m X_{t(s)}
\bigr) \vee\Bigl( \sup_{t(s) \le r \le t(s+h)} \bigl(X_r \m X_{t(s)}\bigr)
\Bigr)\Bigr]\nonumber\\[-8pt]\\[-8pt]
&&\hspace*{159.1pt}{} - \bigl( X_{t(s+h)} \m X_{t(s)} \bigr) \Bigr) \nonumber\\
&=& e^h \Bigl(\Bigl[Z_s^z \vee e^s \Bigl( \sup_{t(s) \le r \le t(s+h)}
\bigl(X_r \m X_{t(s)}\bigr) \Bigr)\Bigr] - e^s \bigl( X_{t(s+h)} \m X_{t(s)}
\bigr) \Bigr) \hspace*{-28pt}\nonumber
\end{eqnarray}
for $s \ge0$ and $h \ge0$ given and fixed. By stationary
independent increments and the scaling property of $X$, we see that
%
%
\begin{eqnarray}
\label{2.17}
\sup_{t(s) \le r \le t(s+h)} \bigl(X_r \m X_{t(s)}\bigr) &=& \sup_{1-e^{-\alpha s} \le r \le1 - e^{-\alpha(s+h)}} (X_r \m
X_{1-e^{-\alpha s}}) \nonumber\\
&\stackrel{\mathrm{law}}{=}& \sup_{0 \le r \le
e^{-\alpha s} (1-e^{-\alpha h})} X_r \nonumber\\[-8pt]\\[-8pt]
& \stackrel{\mathrm{law}}{=}&
\sup_{0 \le r e^{\alpha s} \le1-e^{-\alpha h}} X_{(r
e^{\alpha
s})/e^{\alpha s}} \nonumber\\
&\stackrel{\mathrm{law}}{=}& e^{-s} \sup_{0 \le r
\le
1-e^{-\alpha h}} X_r = e^{-s} S_{t(h)}\nonumber
\end{eqnarray}
and likewise
%
%
\begin{eqnarray}
\label{2.18}
X_{t(s+h)} \m X_{t(s)} &=& X_{1-e^{\alpha(s+h)}} - X_{1-e^{-\alpha
s}} \stackrel{\mathrm{law}}{=} X_{e^{-\alpha s}(1-e^{-\alpha
h})} \nonumber\\[-8pt]\\[-8pt]
&\stackrel{\mathrm{law}}{=}& e^{-s} X_{1-e^{-\alpha h}} =
e^{-s} X_{t(h)}\nonumber
\end{eqnarray}
both being independent from $\cF_{t(s)}^X$. Combining
(\ref{2.16})--(\ref{2.18}), we get
%
%
\begin{equation}
\label{2.19}
\EE\bigl( f ( Z_{s+h}^z ) \vert\cF_{t(s)}^X \bigr) =
\EE\bigl( f \bigl( e^h \bigl(w \vee S_{t(h)} \m X_{t(h)}\bigr) \bigr) \bigr)
\big\vert_{w = Z_s^z}
\end{equation}
for any (bounded) measurable function $f \dvtx \R_+ \rightarrow\R$ from
where all the claims follow by standard means [observe that the
deterministic function on the right-hand side of (\ref{2.19}) does
not depend on $s$ (implying that $Z$ is a time-homogenous Markov
process) as well as that it defines a~continuous and bounded
function of $w$ whenever $f$ is so (Feller property) implying that
$Z$ is a strong Markov process]. This completes the proof.
\end{pf}

Note from (\ref{2.14}) that $Z$ is a transient process (satisfying $Z_s
\rightarrow\infty$ as $s \rightarrow\infty$) having downward jumps only
(since $X$ jumps upward). The state space of $Z$ equals~$\R_+$.

\section{The optimal stopping problem}\label{sec3}

1. From (\ref{2.10}) and (\ref{2.13}), we see that the optimal
prediction problem (\ref{2.3}) reduces to the optimal stopping
problem
%
%
\begin{equation}
\label{3.1}
V = \inf_{0 \le\sigma<\infty} \EE e^{-p \sigma} G(Z_\sigma),
\end{equation}
where the infimum is taken over all stopping times $\sigma$ with
respect to $(\cF_{t(s)}^X)_{s \ge0}$. This optimal stopping problem
is one-dimensional and the horizon is infinite. The exponential term
$(e^{-ps})_{s \ge0}$ in (\ref{3.1}) corresponds to a new (strong)
Markov process $\widetilde Z$ which may be identified with $Z$
killed at rate $p$.

2. To tackle the problem (\ref{3.1}), we need to enable $Z$ to start
at any point in the state space $\R_+$. This can be done using the
result of Proposition \ref{prop1} above, and it leads to the following
variational extension of (\ref{3.1}):
%
%
\begin{equation}
\label{3.2}
V(z) = \inf_{0 \le\sigma<\infty} \EE_z e^{-p \sigma} G(Z_\sigma),
\end{equation}
where the infimum is taken over all stopping times $\sigma$ with
respect to $(\cF_{t(s)}^X)_{s \ge0}$, and the process $Z$ starts at $z$
under $\PP_{ z}$. Moreover, by the result of Proposition \ref{prop1} we
know that $\PP_{ z}$ can be realized by (\ref{2.15}) in terms of $Z^z =
(Z_s^z)_{s \ge0}$ under $\PP$, and this fact will be useful below when
analysing properties of the mapping $z \mapsto V(z)$ on~$\R_+$.

3. Before we turn to a more detailed analysis of the problem
(\ref{3.2}), let us state some basic properties of $G$ and $V$ that
will be useful throughout. Recall that $f(z) \sim g(z)$ as $z
\rightarrow z_0$ means that $\lim_{ z \rightarrow z_0} f(z)/g(z) =
1$ for $z_0 \in[-\infty,\infty]$.
\begin{Proposition}\label{prop2}
The gain function $G$ from
(\ref{2.12}) above and the value function $V$ from (\ref{3.2}) above
satisfy the following properties:
%
%
\begin{eqnarray}
\label{3.3}
&&\hspace*{-6pt}
\begin{tabular}{p{208pt}@{}}
$z \mapsto G(z) \mbox{ is (strictly) increasing and convex on }
\R_+\break \mbox{with } G(0) = \EE S_1^p > 0 ;$
\end{tabular}\\
\label{3.4}
&&z \mapsto V(z) \mbox{ is increasing and continuous on } \R_+
; \\
\label{3.5}
&& z^p \le V(z) \le G(z) \mbox{ for all }
z \in\R_+ ; \\
\label{3.6}
&&G(z) \sim z^p
\mbox{ and } V(z) \sim z^p \mbox{ as } z \rightarrow
\infty.
\end{eqnarray}
\end{Proposition}
\begin{pf}
Equation (\ref{3.3}): recalling that $F_{S_1}$ denotes the
distribution function of~$S_1$, and letting $f_{S_1}$ denote the
density function of $S_1$, we find from the final expression in
(\ref{2.12}) that $G'(z) = p z^{p-1} F_{S_1}(z) > 0$ and $G''(z)
= p (p \m1) \times z^{p-2} F_{S_1}(z) + p z^{p-1} f_{S_1}(z) > 0$
for all $z>0$ implying that $z \mapsto G(z)$ is (strictly)
increasing and convex, respectively. Likewise, we also see from the
middle expression in (\ref{2.12}) that $G(0) = \EE S_1^p > 0$ as
claimed.

Equation (\ref{3.4}): letting $\sigma$ be a given and fixed stopping time, we
see from~(\ref{2.15}) that $z \mapsto Z_\sigma^z$ is increasing so
that $z \mapsto G(Z_\sigma^z)$ is increasing, and the fact that $z
\mapsto V(z)$ is increasing follows directly from the definition
(\ref{3.2}). To show that $z \mapsto V(z)$ is continuous, take $z_1
< z_2$ in $\R_+$ and note by the mean value theorem and (\ref{2.15})
that
%
%
\begin{eqnarray}
\label{3.7}
0 &\le& G(Z_\sigma^{z_2}) \m G(Z_\sigma^{z_1}) = G'(\xi)
(Z_\sigma^{z_2} \m Z_\sigma^{z_1}) \nonumber\\
&=& G'(\xi) e^\sigma\bigl(
z_2 \vee S_{t(\sigma)} \m z_1 \vee S_{t(\sigma)} \bigr)
\\
&\le& p \xi^{p-1} F_{S_1}(\xi)
e^\sigma(z_2 \m z_1),\nonumber
\end{eqnarray}
where $\xi\in(Z_\sigma^{z_1},Z_\sigma^{z_2})$. Since $0 \le\xi
\le e^\sigma(z_2 \vee S_1 - I_1)$, where we set $I_1 = \inf_{ 0
\le t \le1} X_t$, it follows from (\ref{3.7}) that
%
%
\begin{equation}
\label{3.8}\qquad
0 \le\EE e^{-p \sigma} G(Z_\sigma^{z_2}) \m\EE e^{-p
\sigma} G(Z_\sigma^{z_1}) \le p \EE(z_2 \vee S_1 \m I_1)^{p-1}
(z_2 \m z_1) .
\end{equation}
Taking the infimum over all stopping times $\sigma$ it follows that
%
%
\begin{equation}
\label{3.9}
0 \le V(z_2) \m V(z_1) \le K (z_2 \m z_1),
\end{equation}
where $K = p \EE(z_2 \vee S_1 \m I_1)^{p-1} < \infty$. This
implies that $V$ is continuous on~$\R_+$ (as well as Lipschitz
continuous on compact sets in $\R_+$).

Equation (\ref{3.5}): the second inequality is obvious so let us derive the
first inequality. For this, fix any $z \in\R_+$ and note that $G(z)
\ge z^p$ and Jensen's inequality imply that
%
%
\begin{eqnarray}
\label{3.10}
V(z) &\ge&\inf_{0 \le\sigma< \infty} \EE e^{-p \sigma}
(Z_\sigma^z)^p \ge\Bigl( \inf_{0 \le\sigma< \infty} \EE
e^{-\sigma} Z_\sigma^z \Bigr)^{ p} \nonumber\\
&=& \Bigl( \inf_{0 \le\sigma
< \infty} \EE\bigl( z \vee S_{t(\sigma)} \m X_{t(\sigma)}
\bigr) \Bigr)^{ p} \\
&=& \Bigl( \inf_{0 \le\tau\le1}
\EE( z \vee S_\tau\m X_\tau) \Bigr)^{ p}
= z^p\nonumber
\end{eqnarray}
upon using that there is a one-to-one correspondence between
$\sigma$ and $\tau$ as stated following (\ref{2.11}) above. Note
also that for the final equality we use the fact that $\EE X_\tau=
0$ since $X$ is a martingale. This establishes the first inequality
in (\ref{3.5}) as claimed.

Equation (\ref{3.6}): note that (\ref{2.12}) above implies that $G(z)/z^p
\rightarrow1$ as \mbox{$z \rightarrow\infty$}, so that $V(z)/z^p
\rightarrow1$ as $z \rightarrow\infty$ follows by (\ref{3.5}).
This completes the proof.
\end{pf}

4. \textit{Existence of an optimal stopping time}. General theory of
optimal stopping for Markov processes (see, e.g., \cite{PS}) can be
used to establish the existence of an optimal stopping time in the
problem (\ref{3.2}). For this, let $C = \{ z \in\R_+ \dvtx V(z) <
G(z) \}$ denote the (open) continuation set, let $D = \{ z \in
\R_+ \dvtx V(z) = G(z) \}$ denote the (closed) stopping set, and note
that
%
%
\begin{equation}
\label{3.11}
\EE\Bigl( \sup_{s \ge0} e^{-ps} G(Z_s^z) \Bigr) < \infty
\end{equation}
since\vspace*{1pt} $e^{-ps} G(Z_s^z) = e^{-ps} ( (Z_s^z)^p +
\int_{(Z_s^z)^p}^\infty\PP(S_1^p > w) \,dw ) \le(z \vee
S_1 \m I_1)^p + \EE S_1^p$ for all $s \ge0$, and the latter
random variable clearly is integrable for each $z \in\R_+$.
Moreover, by (\ref{3.3}) and (\ref{3.4}) we know that the gain
function $z \mapsto G(z)$ is lower semicontinuous on $\R_+$ and the
value function $z \mapsto V(z)$ is upper semicontinuous on $\R_+$.
Hence, by Corollary 2.9 and Remark 2.10 in \cite{PS}, pages 46--48, we
can conclude that the first entry time of $Z$ into $D$ given by
%
%
\begin{equation}
\label{3.12}
\sigma_{ D} = \inf\{ s \ge0 \dvtx Z_s \in D \}
\end{equation}
is an optimal stopping time in (\ref{3.2}). This stopping time is
not necessarily finite valued [when the set in~(\ref{3.12}) is
empty] and the value $e^{-p \sigma_{ D}} G(Z_{\sigma_{ D}}^z)$
in~(\ref{3.2}) can be formally assigned as $(z \vee S_1 \m X_1)^p$
when $\sigma_{ D} = \infty$ since by~(\ref{2.12}) and (\ref{2.15})
we have
%
%
\begin{equation}
\label{3.13}
e^{-ps} G(Z_s^z) \rightarrow(z \vee S_1 \m X_1)^p
\end{equation}
as $s \rightarrow\infty$. This is in agreement with the usual
hypothesis from general theory introduced to cover the case of
infinite-valued stopping times.

5. In addition to these general facts, it may be noted that the
optimal stopping problem (\ref{3.2}) plays an auxiliary role in
tackling the optimal prediction problem (\ref{2.3}), and it is clear
from our considerations above that we only need to compute $V(z)$
for $z=0$. Thus, if we set $z_* = \inf D$ then either $z_* < \infty$
when $D \ne\varnothing$ (so that $z_* \in D$ since $D$ is closed) or
$z_* = \infty$ when $D = \varnothing$. In the first case (when $D \ne
\varnothing$), the first entry time of $Z$ to $z_*$ given by
%
%
\begin{equation}
\label{3.14}
\sigma_{ z_*} = \inf\{ s \ge0 \dvtx Z_s = z_* \}
\end{equation}
is optimal in (\ref{3.2}) under $\PP_{ z}$ for $z=0$. It should be
recalled here that $Z$ jumps downward only and creeps upward in
$\R_+$ so that $Z$ will hit any point in $(0,\infty)$ with
probability one due to its transience to $+\infty$. Recalling
further the time change (\ref{2.11}) we see that (\ref{3.14})
translates into the fact that the stopping time
%
%
\begin{equation}
\label{3.15}
\tau_* = \inf\{ t \in[0,1] \dvtx S_t \m X_t \ge z_* (1 \m
t)^{1/\alpha} \}
\end{equation}
is optimal in (\ref{2.3}) with $T=1$. In the second case (when $D =
\varnothing$), we see that the optimal stopping time $\sigma_{ z_*}$
in (\ref{3.2}) equals $+\infty$ under $\PP_{ z}$ for $z=0$. In this
case, we have
%
%
\begin{equation}
\label{3.16}
V(z) = \EE(z \vee S_1 \m X_1)^p
\end{equation}
for all $z \in\R_+$ and the time change (\ref{2.11}) implies that
$\tau_* \equiv1$ is optimal in~(\ref{2.3}) with
%
%
\begin{equation}
\label{3.17}
V = \EE(S_1 \m X_1)^p .
\end{equation}
A central question therefore becomes to examine when $[0,z_*)
\subseteq C$ with \mbox{$z_* \in D$} (it will be shown in Section \ref{sec5} below
that $z_*$ cannot be zero). We will tackle this question by forming
a free-boundary problem on $[0,z_*)$ for~$V$ defined in (\ref{3.2}).
For this, we first need to determine the infinitesimal
characteristics of~$Z$.

\section{The free-boundary problem}\label{sec4}

1. The following proposition determines the action of the
infinitesimal generator of the process $Z$ defined in (\ref{2.14})
in terms of the action of the infinitesimal generator of the
reflected process $Y = S \m X$. Below we let $C_b^2(\R_+)$ denote
the class of twice continuously differentiable functions $F \dvtx \R_+
\rightarrow\R$ such that $F'$ and $F''$ are bounded on~$\R_+$.

\begin{Proposition}\label{prop3}
The infinitesimal generator $\LL_Z$ of
the process $Z$ is given by
%
%
\begin{equation}
\label{4.1}
\LL_Z F(z) = z F'(z) + \alpha\LL_Y F(z)
\end{equation}
for any $F \in C_b^2(\R_+)$ satisfying (\ref{4.6}) below,
where $\LL_Y$ denotes the infinitesimal generator of the process
$Y$.
\end{Proposition}
\begin{pf}
By the mean value theorem, we have
%
%
\begin{eqnarray}
\label{4.2}\qquad
\LL_Z F(z) &=& \lim_{s \downarrow0} \frac{1}{s} \EE
\bigl( F(Z_s^z) \m F(z) \bigr) \nonumber\\
&=& \lim_{s \downarrow
0} \frac{1}{s} \EE\bigl( F\bigl(e^s\bigl(z \vee S_{t(s)} \m
X_{t(s)}\bigr) \bigr) - F\bigl(z \vee S_{t(s)} \m X_{t(s)} \bigr)
\nonumber\\
&&\hspace*{106pt}{}
+ F\bigl(z \vee S_{t(s)} \m X_{t(s)}
\bigr) \m F(z) \bigr) \\
&=& \lim_{s \downarrow0}
\frac{e^s \m1}{s} \EE\bigl( F'(\xi_s) \bigl(z \vee S_{t(s)}
\m X_{t(s)}\bigr) \bigr) \nonumber\\
&&{}
+ \lim_{s \downarrow0}
\frac{t(s)}{s} \biggl( \frac{1}{t(s)} \bigl[\EE F\bigl(z \vee
S_{t(s)} \m X_{t(s)} \bigr) \m F(z)\bigr] \biggr) \nonumber\\
&=& z F'(z) + \alpha\LL_Y F(z),\nonumber
\end{eqnarray}
where for the second last limit we use that $(e^s \m1)/s
\rightarrow1$ and $F'(\xi_s) \rightarrow F'(z)$ as $s \downarrow0$
since $\xi_s \in(z \vee S_{t(s)} \m X_{t(s)},e^s (z \vee
S_{t(s)} \m X_{t(s)}))$, and for the last limit we use that $t(s)/s
\rightarrow\alpha$ as $s \downarrow0$ and the result of
Proposition \ref{prop4} below. This completes the proof.
\end{pf}

2. The following proposition determines the action of the
infinitesimal generator of the reflected process $Y = S \m X$. We
refer to the \hyperref[app]{Appendix} for the analogous result in the case of a
general (strictly) stable L\'evy\vadjust{\goodbreak} process~$X$.

\begin{Proposition}\label{prop4}
The infinitesimal generator $\LL_Y$ of
the reflected process $Y = S \m X$ takes any of the following three
forms for $y>0$ given and fixed:
%
%
\begin{eqnarray}
\label{4.3}
&&\mbox{It\^o's form} \nonumber\\
&&\qquad\LL_Y F(y) = \int_0^y
\bigl( F(y \m x) \m F(y) \p F'(y) x \bigr) \frac{c}{x^{1+\alpha}}
\,dx \\
&&\hphantom{\qquad\LL_Y F(y) =}{}+ \frac{c (F(0) \m F(y))}{\alpha y^\alpha} + \frac{c
F'(y)}{(\alpha\m1) y^{\alpha-1}},\nonumber
\end{eqnarray}
\begin{itemize}
\item[]{Riemann--Liouville's form}
\begin{equation}\label{4.4}
\LL_Y F(y)
= \frac{c}{\alpha(\alpha\m1)} \,\frac{d^2}{dy^2} \int_0^y
\frac{F(x)}{(y \m x)^{\alpha-1}} \,dx + \frac{c F(0)}{\alpha
y^\alpha},
\end{equation}
\item[]\mbox{Caputo's form}
 \begin{equation}
\label{4.5} \LL_Y F(y) = \frac{c}{\alpha(\alpha
\m1)} \int_0^y \frac{F''(x)}{(y \m x)^{\alpha- 1}} \,dx,
\end{equation}
\end{itemize}
whenever $F \in C_b^2(\R_+)$ satisfies
%
%
\begin{equation}
\label{4.6}
F'(0+) = 0 \qquad\mbox{(normal reflection)}.
\end{equation}
\end{Proposition}
\begin{pf}
It is enough to establish (\ref{4.3}) since
(\ref{4.4}) and (\ref{4.5}) can then be derived by (repeated)
integration by parts using (\ref{4.6}) (note that the equivalence
of (\ref{4.3})--(\ref{4.5}) under (\ref{4.6}) remain valid for any $F
\in C^1[0,\infty) \cap C^2(0,\infty)$ satisfying $\vert F''(x) \vert
= O(x^{\alpha-2})$ as $x \downarrow0$ since $\alpha\m2 > -1$).
For this, fix $t>0$ and note that by It\^o's formula we have
%
%
\begin{eqnarray}
\label{4.7}
F(Y_t) &=& F(Y_0) + \int_0^t F'(Y_{s-}) \,dY_s \nonumber\\[-8pt]\\[-8pt]
&&{}+ \sum_{0 < s \le t}
\bigl(F(Y_s) \m F(Y_{s-}) \m F'(Y_{s-}) \Delta Y_s \bigr)\nonumber
\end{eqnarray}
since $[Y,Y]^c \equiv0$. Indeed, the latter equality follows by
recalling that $X$ is a quadratic pure jump semimartingale
(i.e., $[X,X]^c=0$) since it is a L\'evy process with no
Brownian component (see \cite{Pr}, page 71), the process $S$ is a
quadratic pure jump semimartingale since it is of bounded variation
(see Theorem 26 in \cite{Pr}, page 71), and the sum/difference of two
quadratic pure jump semimartingales is a quadratic pure jump
semimartingale (this can be easily verified using Theorem 28 in
\cite{Pr}, page 75, e.g.).

Since $X$ jumps upward and creeps downward, it follows that $dS_s
= \Delta S_s$ in terms of a suggestive notation, and hence from
(\ref{4.7}) we get
%
%
\begin{equation}
\label{4.8}
F(Y_t)\,{=}\,F(Y_0)\,{+}\,M_t\,{+}\,\sum_{0 < s \le t}
\bigl(F(Y_{s-}\,{+}\,\Delta Y_s)\,{-}\,F(Y_{s-})\,{+}\,F'(Y_{s-}) \Delta X_s \bigr),\hspace*{-40pt}
\end{equation}
where $M_t = - \int_0^t F'(Y_{s-}) \,dX_s$ is a local martingale for
$t \ge0$. By the BDG inequality (see, e.g., \cite{PS}, page 63)
combined with the facts that $F'$ is bounded on $\R_+$ and $\EE
[X,X]_t^q < \infty$ with $q=1/2$ since $[X,X]$ is a stable process
of index $\alpha/2>q$ [with L\'evy measure $c \,dx/(2
x^{1+\alpha/2})$ as is easily verified directly from definition] it
follows that ${\EE\sup_{ 0 \le s \le t} }\vert M_s \vert< \infty$
and hence $M$ is a~martingale. The right-hand side of this identity
can be further rewritten as follows
%
%
\begin{eqnarray}
\label{4.9}\qquad
F(Y_t) &=& F(Y_0) + M_t \nonumber\\
&&{}+ \sum_{0 < s \le t} \bigl( [ F(Y_{s-}
- \Delta X_s)  \nonumber\\[-8pt]\\[-8pt]
&&\hspace*{42.2pt}{}\m F(Y_{s-})\p F'(Y_{s-}) \Delta X_s ]
I(\Delta X_s \le Y_{s-}) \nonumber\\
&&\hspace*{39.3pt}{}
+ [F(0)
\m F(Y_{s-}) + F'(Y_{s-}) \Delta X_s ] I(\Delta X_s
> Y_{s-}) \bigr)\nonumber
\end{eqnarray}
upon using that $\Delta X_s \le Y_{s-}$ if and only if $X_s \le
S_{s-}$ so that $\Delta S_s=0$, and $\Delta X_s > Y_{s-}$ if and
only if $X_s>S_{s-}$ so that $S_s = X_s$, that is, $Y_s=0$. Taking~%
$\EE_y$ on both sides of (\ref{4.9}), where $\PP_{ y}$ denotes a
probability measure under which $Y_0=y$, and applying the
compensation formula (see, e.g., \cite{RY}, page~475) we find that
%
%
\begin{eqnarray}
\label{4.10}
&&\EE_y F(Y_t) \m F(y) \nonumber\\
&&\qquad= \EE_y \biggl[\int_0^t ds \biggl(
\int_0^{Y_s} [ F(Y_s \m x) \m F(Y_s) \p F'(Y_s)
x ] \nu(dx)\\
&&\qquad\quad\hspace*{61.5pt}{}
+ \int_{Y_s}^\infty[
F(0) \m F(Y_s) \p F'(Y_s) x ] \nu(dx) \biggr) \biggr]\nonumber
\end{eqnarray}
for all $y>0$. The applicability of this formula (see, e.g.,
\cite{Ky}, page 97) follows from the facts that $\vert F'(y)\vert\,{\le}\,C y$ and
$\vert F''(y)\vert\,{\le}\,C$ for all $y\,{\ge}\,0$ with some \mbox{$C\,{>}\,0$}
so that the mean value theorem yields the existence of
$\xi_{s,x} \in(Y_s \m x,Y_s)$ and $\eta_s \in(0,Y_s)$ such that
%
%
\begin{eqnarray}
\label{4.11}
&&\EE_y \biggl[\int_0^t ds \biggl(\int_0^{Y_s} \vert F(Y_s \m
x) \m F(Y_s) \p F'(Y_s) x \vert\nu(dx) \nonumber\\
&&\quad\hspace*{49.5pt}{}
+ \int_{Y_s}^\infty\vert F(0) \m F(Y_s)
\p F'(Y_s) x \vert\nu(dx) \biggr) \biggr] \nonumber\\
&&\qquad\le\EE_y \biggl[ \int_0^t ds \biggl(\int_0^{Y_s} \frac{1}
{2} \vert F''(\xi_{s,x}) \vert x^2 \frac{c}{x^{1+\alpha}}
\,dx \nonumber\\[-8pt]\\[-8pt]
&&\qquad\hspace*{59.2pt}{}
+ \int_{Y_s}^\infty\bigl(\vert F'(\eta_s)
\vert Y_s + \vert F'(Y_s) \vert x \bigr) \frac{c}{x^{1+
\alpha}} \,dx \biggr) \biggr] \nonumber\\
&&\qquad\le c \EE_y \biggl[\int_0^t
ds \biggl( \frac{C}{2(2-\alpha)} + \frac{C}{\alpha} + \frac{C}
{\alpha-1} \biggr) Y_s^{2-\alpha} \biggr] \nonumber\\
&&\qquad\le c \biggl(
\frac{C}{2(2-\alpha)} + \frac{C}{\alpha} + \frac{C} {\alpha-1}
\biggr) \frac{\alpha}{2} t^{2/\alpha} \EE_y(S_1 \m I_1)^{2-
\alpha} < \infty\nonumber
\end{eqnarray}
since $2 \m\alpha\in(0,\alpha)$ and where we also use the scaling
property of $X$. Dividing both sides of (\ref{4.10}) by $t$, letting
$t \downarrow0$ and using the dominated\vadjust{\goodbreak} convergence theorem, we get
%
%
\begin{eqnarray}
\label{4.12}
\LL_Y F(y) &=& \int_0^y [ F(y \m x) \m F(y) \p F'(y) x
] \nu(dx) \nonumber\\[-8pt]\\[-8pt]
&&{}
+ [ F(0) \m F(y) ]
\int_y^\infty\nu(dx) + F'(y) \int_y^\infty x \nu(dx),\nonumber
\end{eqnarray}
which is easily verified to be equal to the right-hand side of
(\ref{4.3}) for all $y>0$ upon using (\ref{2.4}). This completes the
proof.
\end{pf}

3. It will be shown in Section \ref{sec5} below that the continuation set
$C$ in the optimal stopping problem (\ref{3.2}) always contains the
interval $[0,\eps)$ for some $\eps> 0$ sufficiently small, so that
the optimal stopping point $z_*$ from (\ref{3.14}) is always strictly
larger than zero. Moreover, we now show that the value function $V$
from (\ref{3.2}) is smooth from the left at $z_*$ whenever $D \ne
\varnothing$.

\begin{Proposition}[(Smooth fit)]\label{prop5}
If the optimal stopping
point $z_*$ from~(\ref{3.14}) is finite, then the value function $V$
from~(\ref{3.2}) is differentiable from the left at $z_*$ and we
have
%
%
\begin{equation}
\label{4.13}
V_-'(z_*) = G'(z_*) .
\end{equation}
\end{Proposition}
\begin{pf}
To simplify the notation, let us write $b$ in place
of $z_*$. Then $[0,b) \subseteq C$ and $b \in D$ so that
$V(b)=G(b)$. Hence, $(V(b \m\eps) \m V(b))/(-\eps) \ge(G(b \m\eps)
\m G(b))/(-\eps)$ for all $\eps>0$ sufficiently small, and letting
$\eps\downarrow0$ we obtain
%
%
\begin{equation}
\label{4.14}
\liminf_{\eps\downarrow0} \frac{V(b \m\eps) \m V(b)}{-\eps}
\ge G'(b) .
\end{equation}

To derive a reverse inequality, note that the stopping time
%
%
\begin{equation}
\label{4.15}
\sigma_\eps= \inf\{ s \ge0 \dvtx Z_s^{b-\eps} \ge b \}
\end{equation}
is optimal for $V(b \m\eps)$ under $\PP$ (recall that $Z$ creeps
upward). Hence, by the mean value theorem we find that
%
%
\begin{eqnarray}
\label{4.16}
&&
V(b \m\eps) \m V(b) \nonumber\\
&&\qquad\ge\EE( e^{-p \sigma_\eps}
G(Z_{\sigma_\eps}^{b-\eps}) ) - \EE(e^{-p
\sigma_\eps} G(Z_{\sigma_\eps}^b) ) \nonumber\\
&&\qquad= \EE\bigl(
e^{-p \sigma_\eps} G'(\xi_\eps) (Z_{\sigma_\eps}^{b-\eps}
\m Z_{\sigma_\eps}^b) \bigr) \\
&&\qquad= \EE\bigl(
e^{-p \sigma_\eps} G'(\xi_\eps) \bigl( e^{\sigma_\eps}
\bigl( (b \m\eps) \vee S_{t(\sigma_\eps)} - b \vee
S_{t(\sigma_\eps)} \bigr) \bigr) \bigr) \nonumber\\
&&\qquad\ge
-\eps\EE\bigl( e^{-p \sigma_\eps} G'(\xi_\eps)
e^{\sigma_\eps} I\bigl(S_{t(\sigma_\eps)} < b \bigr) \bigr),\nonumber
\end{eqnarray}
where $\xi_\eps\in(Z_{\sigma_\eps}^{b-\eps}, Z_{\sigma_\eps}^b)$
for $\eps\in(0,b)$.

We claim that $\sigma_\eps\rightarrow0$ $\PP\mbox{-a.s.}$ as $\eps
\downarrow0$. Indeed, setting
%
%
\begin{eqnarray}
\label{4.17}
\rho_\eps&=& \inf\bigl\{ s \ge0 \dvtx (b \m\eps) \vee S_{t(s)} \m
X_{t(s)} \ge b \bigr\} ,\\
\tau_\eps&=& \inf\{ t \ge0 \dvtx
(b \m\eps) \vee S_t \m X_t \ge b \}
\end{eqnarray}
we see that $\sigma_\eps\le\rho_\eps$ and $\rho_\eps=
t^{-1}(\tau_\eps)$ for all $\eps> 0$. Since $t^{-1}(0+) = 0$ it is
therefore sufficient to show that $\tau_\eps\rightarrow0$
$\PP\mbox{-a.s.}$ as $\eps\downarrow0$. For this, note that
%
%
\begin{equation}
\label{4.19}
\tau_\eps\le\inf\{ t \ge0 \dvtx (b \m\eps) \m X_t \ge b \}
= \inf\{ t \ge0 \dvtx X_t \le-\eps\} =: \gamma_\eps
\end{equation}
and $\gamma_\eps\downarrow0$ $\PP\mbox{-a.s.}$ as $\eps\downarrow
0$ since the starting point $0$ of $X$ is regular for $(-\infty,0)$.
Hence, $\sigma_\eps\rightarrow0$ $\PP\mbox{-a.s.}$ for $\eps
\downarrow0$ as claimed.

Dividing both sides of (\ref{4.16}) by $-\eps$, letting $\eps
\downarrow0$, and using the dominated convergence theorem [upon
noting that $\xi_\eps\le b \p(Z_{\sigma_\eps}^b \m
Z_{\sigma_\eps}^{b-\eps}) \le b \p\eps e^{\sigma_\eps} \le(b \p
\eps) e^{\sigma_\eps}$ and recalling that $G'(z) = p z^{p-1}
F_{S_1}(z) \le2 z^{p-1}$ for all $z \ge0$ so that $0 \le e^{-p
\sigma_\eps} G'(\xi_\eps) e^{\sigma_\eps} I(S_{t(\sigma_\eps)}
< b) \le2 e^{(-p+1)\sigma_\eps}$ $(b \p\eps)^{p-1}
e^{(p-1)\sigma_\eps} = 2 (b \p\eps)^{p-1} \le2 (b \p1)^{p-1}$ as
$\eps\downarrow0$], we get
%
%
\begin{equation}
\label{4.20}
\limsup_{\eps\downarrow0} \frac{V(b \m\eps) \m V(b)}{-\eps}
\le G'(b) .
\end{equation}
Combining (\ref{4.14}) and (\ref{4.20}), we see that $V$ is
differentiable from the left at~$b$ and that (\ref{4.13}) holds as
claimed. This completes the proof.
\end{pf}

4. Returning to the case when $[0,z_*) \subseteq C$ with $z_* \in
D$, recalling the general fact on the killed Dirichlet problem
(which suggests that $z \mapsto V(z) = \EE_z e^{-p
\sigma_{ z_*}} G(Z_{\sigma_{z_*}})$ should solve $\LL_Z V =
p V$ in $[0,z_*)$ due to the strong Markov property of $Z$; see,
e.g., \cite{PS}, pages 130--132), and making use of the facts from
Propositions \ref{prop3}--\ref{prop5}, we can formulate the following free-boundary
problem for the value function $V$ defined in (\ref{3.2}) above:
%
%
\begin{eqnarray}
\label{4.21}
&&z V'(z) + \alpha\LL_Y V(z) - p V(z) = 0 \qquad\mbox{for } z
\in[0,z_*), \\
\label{4.22}
&&V(z_*) = G(z_*) \qquad\mbox{(instantaneous stopping)}, \\
\label{4.23}
&&V'(z_*) = G'(z_*)\qquad
\mbox{(smooth fit)}, \\
\label{4.24}
&&V'(0) = 0\qquad
\mbox{(normal reflection)},
\end{eqnarray}
where $z_* \in(0,\infty)$ is the (unknown) boundary point to be
found along with~$V$ on $[0,z_*)$. Whilst the infinitesimal
generator $\LL_Y$ in (\ref{4.21}) can take any of the three forms
(\ref{4.3})--(\ref{4.5}) from Proposition \ref{prop4}, it turns out that the
Caputo form (\ref{4.5}) is most convenient for the analysis of the
problem (\ref{4.21})--(\ref{4.24}) to be performed.

For this reason, let us rewrite (\ref{4.21}) in the
Caputo form as
%
%
\begin{equation}
\label{4.25}
z F'(z) + \frac{c}{\alpha\m1} \int_0^z \frac{F''(x)}{(z \m
x)^{\alpha-1}} \,dx - p F(z) = 0
\end{equation}
for $z \in(0,b]$ and $F \dvtx [0,b] \rightarrow\R$ with $b \in
(0,\infty)$ given and fixed. The proof of Proposition \ref{prop6} below shows
that the natural solution space for this equation is one-dimensional
[once $F'(0)$ is set to $0$]. More precisely, let $S_b$ denote the
class of functions $F \dvtx [0,b] \rightarrow\R$ satisfying the
following three conditions:
%
%
\begin{eqnarray}
\label{4.26}
&&F \in C^1[0,b] \cap C^2(0,b], \\
\label{4.27}
&&\vert F''(z)\vert = O(z^{\alpha-2}) \qquad\mbox{as } z \downarrow0, \\
\label{4.28}
&&F'(0) = 0 .
\end{eqnarray}
Note that $F''$ is assumed to exist (and be continuous) on $(0,b]$
but may be unbounded (locally at zero). Note also that
(\ref{4.26})--(\ref{4.28}) imply that $\vert F'(z) \vert=
O(z^{\alpha-1})$ as $z \downarrow0$. For further reference, let us
also recall the following well-known identity (see, e.g., (3.191) in
\cite{GR}, page 333, and (6.2.2) in \cite{AS}, page 258):
%
%
\begin{equation}
\label{4.29}
\int_0^z x^{\mu-1} (z \m x)^{\nu-1} \,dx = z^{\mu+\nu-1}
\frac{\Gamma(\mu) \Gamma(\nu)}{\Gamma(\mu\p\nu)}
\end{equation}
for $\mu>0$ and $\nu>0$.
\begin{Proposition}\label{prop6}
The equation (\ref{4.25}) has a unique
solution $F$ in $S_b$ satisfying
%
%
\begin{equation}
\label{4.30}
F(0) = a_0,
\end{equation}
whenever $a_0 \in\R$ is given and fixed. Moreover, the
following explicit representation is valid:
%
%
\begin{equation}
\label{4.31}
F(z) = a_0 \sum_{n=0}^\infty\frac{1}{(-c \Gamma(-\alpha))^n}
\biggl( - \frac{p}{\alpha}\biggr)_{ n} \frac{z^{\alpha n}}
{\Gamma(\alpha n \p1)}
\end{equation}
for $z \in[0,b]$ where $(q)_n = q (q \p1) \cdots(q \p n
\m1)$ for $n \ge1$ and $(q)_0 = 1$ with $q = -p/\alpha$.
\end{Proposition}
\begin{pf}
1. \textit{Uniqueness.} We will establish the
uniqueness of solution by reducing the integro-differential equation
(\ref{4.25}) to a Volterra integral equation of the second kind. For
this, let us introduce the following substitution in
(\ref{4.25}):
%
%
\begin{equation}
\label{4.32}
\varphi(z) = \int_0^z \frac{F''(x)}{(z \m x)^{\alpha-1}} \,dx
\end{equation}
for $z>0$ upon extending $F$ from $[0,b]$ to a bounded $C^2$
function on $(0,\infty)$ with bounded support in $\R_+$. Let
$\LL[f](\lambda) = \int_0^\infty e^{-\lambda x} f(x) \,dx$ denote
the Laplace transform of a function $f \dvtx \R_+ \rightarrow\R$ with
$\lambda>0$, and let $\LL^{-1}$ denote the inverse Laplace
transform. By (\ref{4.27}) and (\ref{4.29}) we see that
$\LL[\varphi](\lambda)$ is well defined and finite for all
$\lambda>0$. Applying first $\LL$ and then $\LL^{-1}$ on both sides
of (\ref{4.32}) using the well-known properties (i) $\LL[\int_0^x
f_1(y) f_2(x \m y) \,dy](\lambda) = \LL[f_1](\lambda)
\LL[f_2](\lambda)$, (ii)~$\LL[f''](\lambda) = \lambda^2
\LL[f](\lambda) \m\lambda f(0) \m f'(0)$ and (iii)
$\LL[x^\rho](\lambda) = \Gamma(\rho\p1)/\lambda^{\rho+1}$ for
$\rho> -1$, one finds using (\ref{4.28}) that
%
%
\begin{equation}
\label{4.33}
F(z) = \frac{1}{\Gamma(\alpha) \Gamma(2 \m\alpha)} \int_0^z
(z \m x)^{\alpha-1} \varphi(x) \,dx + F(0)
\end{equation}
for $z \in(0,b]$. Inserting this expression back into (\ref{4.25}),
we obtain
%
%
\begin{equation}
\label{4.34}
\int_0^z K(z,x) \varphi(x) \,dx + \varphi(z) = \psi,
\end{equation}
where $K$ and $\psi$ are given by
%
%
\begin{eqnarray}
\label{4.35}
K(z,x) &=& \frac{\alpha\m1}{c \Gamma(\alpha) \Gamma(2
\m\alpha)} \frac{(\alpha\m1\m p) z \p p x}{(z \m
x)^{2-\alpha}}, \\
\label{4.36}
\psi &=& \frac{p
(\alpha\m1)}{c} F(0)
\end{eqnarray}
for $z \in(0,b]$ and $x \in(0,z)$. We may now recognise
(\ref{4.34}) as a Volterra integral equation of the second kind with
a weakly singular kernel $K$ [the kernel is said to be weakly
singular since the exponent $2 \m\alpha$ in the singular term $(z
\m x)^{2-\alpha}$ belongs to the interval $(0,1)$]. Moreover, since
$\psi$ defines a~bounded function on $[0,b]$, it is well known (see,
e.g., \cite{Ho}, Theorem 7, page~35) that the equation (\ref{4.34})
can have at most one solution $\varphi$ (in the class of locally
integrable functions), and by means of the identity (\ref{4.33})
this fact translates directly into the uniqueness of solution for
(\ref{4.25}) as claimed. This completes the first part
of the proof.

2. \textit{Existence.} Seeking a solution to (\ref{4.25}) of the form
%
%
\begin{equation}
\label{4.37}
F(z) = \sum_{n=0}^\infty a_n z^{\beta n + \gamma}
\end{equation}
and inserting it into (\ref{4.25}) upon differentiating and
integrating formally term by term and making use of the identity
(\ref{4.29}), a lengthy but straightforward calculation shows that
$\beta= \alpha$, $\gamma= 0$ and the series coefficients satisfy
%
%
\begin{equation}
\label{4.38}
a_{n+1} = \frac{1}{c \Gamma(-\alpha)} \biggl( \frac{p}{\alpha}
\m n \biggr) \frac{\Gamma(\alpha n \p1)}{\Gamma(\alpha(n \p1)
\p1)} a_n
\end{equation}
for $n = 0, 1, \ldots.$ This yields the candidate series
representation (\ref{4.31}). Moreover, setting $b_n = (1/(-c
\Gamma( -\alpha))^n) (-p/\alpha)_n (z^{\alpha n}
/\Gamma(\alpha n \p1))$ for $n \ge1$ and using the well-known fact
that $\Gamma(\alpha n \p1) /\Gamma(\alpha(n \p1) \p1) \sim
(\alpha n)^{-\alpha}$ as $n \rightarrow\infty$ (see (6.1.47) in
\cite{AS}, page 257), it is easily verified\vadjust{\goodbreak} that $b_{n+1}/b_n
\rightarrow0$ as $n \rightarrow\infty$. Hence, by the ratio test we
can conclude that the series in (\ref{4.31}) converges absolutely
for every $z \in[0,b]$. A direct verification also shows that the
function $F$ defined by the series in (\ref{4.31}) belongs to $S_b$.
These facts justify the formal steps leading to (\ref{4.38}) above,
and the proof is complete.
\end{pf}

5. Before we continue our analysis of the free-boundary problem
(\ref{4.21})--(\ref{4.24}), let us make precise the following
consequence of It\^o's formula and the optional sampling theorem.
Note that $G$ satisfies both (\ref{4.39}) and (\ref{4.40}) below
since $\vert G''(z) \vert= O(z^{p+\alpha-3})$ as $z \downarrow0$
and $\vert G''(z) \vert= O(z^{p-2})$ as $z \uparrow\infty$. This
is easily seen upon recalling the expression for $G''$ from the
proof of~(\ref{3.3}) above and using the asymptotic relations
(\ref{5.14}), (\ref{5.15}) and (\ref{5.17}) below. Recall also that
$F$ from Proposition \ref{prop6} satisfies (\ref{4.39}) below.
\begin{Proposition}\label{prop7}
Let $F \dvtx \R_+ \rightarrow\R$ be a
function from $C^1[0,\infty) \cap C^2(0,\infty)$ satisfying
%
%
\begin{equation}
\label{4.39}
\vert F''(z) \vert= O(z^{\alpha-2}) \qquad\mbox{as } z
\downarrow0 \quad\mbox{and}\quad F'(0)=0 .
\end{equation}
Let $\sigma$ be a stopping time of $Z$ such that either
%
%
\begin{equation}
\label{4.40}
\vert F''(z) \vert= O(z^\beta) \qquad\mbox{as } z \uparrow
\infty\mbox{ for some } \beta< \alpha\m2
\end{equation}
and $\sigma\le k$ for some $k \ge1$, or $\sigma\le
\sigma_m$ for some $m \ge1$ where $\sigma_m = \inf\{ s \ge0
\dvtx\allowbreak
Z_s = m \}$. Then the following identity holds:
%
%
\begin{equation}
\label{4.41}
\EE_z e^{-p \sigma} F(Z_\sigma) = F(z) + \EE_z \int_0^\sigma
e^{-ps} (\LL_Z F \m p F)(Z_s) \,ds
\end{equation}
for all $z \in\R_+$.
\end{Proposition}
\begin{pf}
Under $\PP_{ z}$ with $z \in\R_+$ by It\^o's
formula, we get
%
%
\begin{eqnarray}
\label{4.42}
e^{-ps} F(Z_s) &=& F(z) - p \int_0^s e^{-pr} F(Z_r) \,dr
\nonumber\\[-8pt]\\[-8pt]
&&{} +\int_0^s e^{-pr} F'(Z_r) Z_r \,dr + M_s + J_s,\nonumber
\end{eqnarray}
where $M_s = - \int_0^s e^{(-p+1)r} F'(Z_{r-}) \,dX_{t(r)}$ is a
local martingale and $J_s =\break \sum_{0 < r \le s} e^{-pr} ( F(Z_r)
\m F(Z_{r-}) \p e^r F'(Z_{r-}) \Delta X_{t(r)} )$ for $s \ge0$
[upon noting that $dZ_r = Z_r \,dr + e^r \,dY_{t(r)}$ from
(\ref{2.15}) above]. Note also that $Z$ is a quadratic pure jump
semimartingale (i.e., $[Z,Z]^c=0$) for the reasons outlined
following (\ref{4.7}) above. Note further that similarly to
(\ref{4.9}) we find that
%
%
\begin{eqnarray}\quad
\label{4.43}
J_s &=& \sum_{0 < r \le s} e^{-pr} \bigl(\bigl[ F\bigl(e^r Y_{t(r)-} -
e^r \Delta X_{t(r)}\bigr) \nonumber\\
&&\hspace*{53.82pt}{}- F\bigl(e^r Y_{t(r)-}\bigr) + e^r
F'\bigl(e^r Y_{t(r)-}\bigr) \Delta X_{t(r)} \bigr]\nonumber\\
&&\hspace*{50.4pt}{}\times I\bigl(\Delta X_{t(r)}
\le Y_{t(r)-}\bigr)\\
&&\hspace*{50.4pt}{}+ \bigl[ F(0) - F\bigl(e^r Y_{t(r)-}\bigr)
+ e^r F'\bigl(e^r Y_{t(r)-}\bigr) \Delta X_{t(r)} \bigr]\nonumber\\
&&\hspace*{171.3pt}{}\times I\bigl(\Delta
X_{t(r)} > Y_{t(r)-}\bigr) \bigr)\nonumber
\end{eqnarray}
upon using that $\Delta X_{t(r)}\,{\le}\,Y_{t(r)-}$ if and only if
$X_{t(r)}\,{\le}\,S_{t(r)-}$ so that \mbox{$\Delta S_{t(r)}\,{=}\,0$}, and $\Delta
X_{t(r)} > Y_{t(r)-}$ if and only if $X_{t(r)} > S_{t(r)-}$ so that
$S_{t(r)} = X_{t(r)}$, that is, $Y_{t(r)}=0$. Setting $v=$ $t(r)$ this
further reads
%
%
\begin{eqnarray}
\label{4.44}\qquad
J_s &=& \sum_{0 < v \le t(s)} e^{-p t^{-1}(v)} \nonumber\\
&&\hspace*{33.6pt}{}\times\bigl(\bigl[ F\bigl(e^{t^{-1}
(v)} Y_{v-} - e^{t^{-1}(v)} \Delta X_v\bigr) \nonumber\\
&&\hspace*{-41pt}\hspace*{93.7pt}{}- F\bigl(e^{t^{-1}(v)}
Y_{v-}\bigr) + e^{t^{-1}(v)} F'\bigl(e^{t^{-1}(v)} Y_{v-}\bigr) \Delta X_v
\bigr]\nonumber\\[-8pt]\\[-8pt]
&&\hspace*{-41pt}\hspace*{89.9pt}{}\times
I(\Delta X_v \le Y_{v-}) \nonumber\\
&&\hspace*{-41pt}\hspace*{89.9pt}{} + \bigl[ F(0) - F\bigl(e^{t^{-1}(v)}
Y_{v-}\bigr) + e^{t^{-1}(v)} F'\bigl(e^{t^{-1}(v)} Y_{v-}\bigr) \Delta X_v
\bigr]\nonumber\\
&&\hspace*{216.4pt}{}\times
I(\Delta X_v > Y_{v-}) \bigr).\nonumber
\end{eqnarray}
The compensator $K$ of $J$ is given by
%
%
\begin{eqnarray}
\label{4.45}
K_s &=& \int_0^{t(s)} e^{-p t^{-1}(v)} \,dv \nonumber\\
&&{}\times\biggl( \int_0^{Y_v}
\bigl[ F\bigl(e^{t^{-1}(v)} Y_v - e^{t^{-1}(v)} x\bigr)
\nonumber\\[-8pt]\\[-8pt]
&&\hspace*{39.9pt}{}-
F\bigl(e^{t^{-1}(v)} Y_v\bigr) + e^{t^{-1}(v)} F'\bigl(e^{t^{-1}(v)} Y_v\bigr)
x \bigr] \nu(dx) \nonumber\\
&&\hspace*{18pt}{}+ \int_{Y_v}^\infty\bigl[ F(0)
- F\bigl(e^{t^{-1}(v)} Y_v\bigr) + e^{t^{-1}(v)} F'\bigl(e^{t^{-1}(v)} Y_v\bigr)
x \bigr] \nu(dx) \biggr).\hspace*{-28pt}\nonumber
\end{eqnarray}
Setting $r=t^{-1}(v)$ and $y = e^r x$ we see that $dv = \alpha
e^{-\alpha r} \,dr$ and $dx = e^{-r} \,dy$ so that $\nu(dx) = c \,dx /
x^{1+\alpha} = (e^{(1+\alpha)r} c \,dx) / y^{1+\alpha} = (e^{\alpha
r} c \,dy) / y^{1+\alpha} = e^{\alpha r} \nu(dy)$. This shows that
%
%
\begin{eqnarray}
\label{4.46}\quad
K_s &=& \alpha\int_0^s e^{-p r} \,dr \biggl(\int_0^{Z_r} [
F(Z_r \m y) - F(Z_r) + F'(Z_r) y ] \nu(dy) \nonumber\\
&&\hspace*{79.3pt}{} +\int_{Z_r}^\infty[ F(0) - F(Z_r) + F'(Z_r) y ]
\nu(dy) \biggr) \\
&=& \alpha\int_0^s e^{-p r} \LL_Y
F(Z_r) \,dr\nonumber
\end{eqnarray}
upon recalling the argument following (\ref{4.12}) above to obtain
the final equality [where $\LL_Y F$ denotes the action of $\LL_Y$ on
$F$ given by the right-hand side of (\ref{4.3})--(\ref{4.5})].

If $m \ge1$ is given and fixed then (\ref{4.39}) implies the
existence of $C>0$ such that $\vert F'(z) \vert\le C
z^{\alpha-1}$ and $\vert F''(z) \vert\le C z^{\alpha-2}$ for
all $z \in(0,m]$. This combined with the mean value theorem yields
the existence of $\xi_{r,y} \in(Z_r \m y,Z_r)$ and $\eta_r \in
(0,Z_r)$ such that
%
%
\begin{eqnarray}
\label{4.47}
&&\EE_z \biggl[ \int_0^{s \wedge\sigma_m} e^{-p r} \,dr
\biggl( \int_0^{Z_r} \vert F(Z_r \m y) - F(Z_r) + F'(Z_r)
y \vert\nu(dy) \nonumber\\
&&\quad\hspace*{92.5pt}{}
+\int_{Z_r}^\infty
\vert F(0) - F(Z_r) + F'(Z_r) y \vert\nu(dy)
\biggr) \biggr] \nonumber\\
&&\qquad\le\EE_z \biggl[ \int_0^{s \wedge\sigma_m}
e^{-p r} \,dr \biggl( \int_0^{Z_r} \frac{1}{2} \vert
F''(\xi_{r,y}) \vert y^2 \frac{c}{y^{1+\alpha}} \,dy \nonumber\\
&&\qquad\hspace*{98pt}{}
+ \int_{Z_r}^\infty\bigl(\vert F'(\eta_r) \vert Z_r +
\vert F'(Z_r) \vert y \bigr) \frac{c}{y^{1+\alpha}} \,dy \biggr)
\biggr] \\
&&\qquad\le c \EE_z \biggl[ \int_0^{s \wedge\sigma_m}
e^{-p r} \,dr \biggl( \frac{C}{2} \int_0^{Z_r} (Z_r
\m y)^{\alpha-2} y^{1-\alpha} \,dy \nonumber\\
&&\qquad\hspace*{103.8pt}{}
+ C
Z_r^\alpha\int_{Z_r}^\infty y^{-1-\alpha} \,dy + C Z_r^{
\alpha-1} \int_{Z_r}^\infty y^{-\alpha} \,dy \biggr) \biggr] \hspace*{-20pt}\nonumber\\
&&\qquad= c \biggl( \frac{C}{2} \Gamma(2 \m\alpha)
\Gamma(\alpha\m1) + \frac{C}{\alpha} + \frac{C}{\alpha-1}
\biggr) \EE_z \biggl[ \int_0^{s \wedge\sigma_m} e^{-p r}
\,dr \biggr] < \infty\nonumber
\end{eqnarray}
upon using (\ref{4.29}) in the final equality. It follows that $N_{s
\wedge\sigma_m} := J_{s \wedge\sigma_m} \m K_{s \wedge\sigma_m}$
is a martingale under $\PP_z$ for $s \ge0$ (see, e.g., \cite{Ky}, page
97). This shows that $N := J \m K$ is a local martingale [with
$(\sigma_m)_{m \ge1}$ as a localization sequence of stopping
times].

Let $\sigma$ be a stopping time of $Z$ such that $\sigma\le
\sigma_m$ for some $m \ge1$. Choose a~localization sequence of
stopping times $(\rho_n)_{n \ge1}$ for the local martingale~$M$.
Subtracting and adding $K_s$ on the right-hand side of (\ref{4.42}),
replacing~$s$ by $\sigma\wedge\rho_n$, taking $\EE_z$ on both
sides and applying the optional sampling theorem, we obtain
%
%
\begin{equation}
\label{4.48}
\EE_z e^{-p ({\sigma\wedge\rho_n})} F(Z_{\sigma\wedge\rho_n})
= F(z) + \EE_z \int_0^{\sigma\wedge\rho_n} e^{-pr} (\LL_Z
F \m p F)(Z_r) \,dr\hspace*{-25pt}
\end{equation}
for all $z \in\R_+$ and all $n \ge1$ [upon recalling (\ref{4.46})
and the action of $\LL_Z$ in~(\ref{4.1}) above]. Moreover, it is
easily seen from (\ref{4.5}) using (\ref{4.39}) and (\ref{4.29})
that $z \mapsto\LL_YF(z)$ is bounded on $[0,m]$ (and so are $F$ and
$F'$ by continuity). Letting $n \rightarrow\infty$ in (\ref{4.48})
and using the dominated convergence theorem we see that (\ref{4.41})
holds as claimed in this case.

Let us now assume that (\ref{4.40}) holds with $\sigma\le k$ for
some $k \ge1$. Choose again a localization sequence of stopping
times $(\rho_n)_{n \ge1}$, however, this time for both the local
martingale $M$ and and the local martingale $N$. Subtracting and
adding $K_s$ on the right-hand side of (\ref{4.42}), replacing $s$
by $\sigma\wedge\rho_n$,~ta\-king~$\EE_z$ on both sides and applying
the optional sampling theorem, we again obtain (\ref{4.48}) for all
$z \in\R_+$ and all $n \ge1$. Moreover, it is easily seen from
(\ref{4.5}) using (\ref{4.39})$+$(\ref{4.40}) and (\ref{4.29}) that
$\vert\LL_YF(z) \vert\le C_3 (1 \p z^{\beta+2-\alpha})$ for all $z
\in\R_+$ with some $C_3>0$. Likewise, it is easily verified that
(\ref{4.39}) and (\ref{4.40}) imply that $\vert F(z) \vert\le C_4
(1 \p z^{\beta+2})$ and $\vert F'(z) \vert\le C_5 (1 \p
z^{\beta+1})$ for all $z \in\R_+$ with some $C_4>0$ and $C_5>0$.
Hence, we see that there exists $C_6>0$ such that
%
%
\begin{eqnarray}
\label{4.49}
&&\vert F(Z_s^z) \vert+ \vert(\LL_Z F \m p F)(Z_s^z) \vert
\nonumber\\
&&\qquad\le C_6 \bigl(1 \p(Z_s^z)^{\beta+2}\bigr) \\
&&\qquad\le C_6 \bigl( 1 + e^{k(\beta
+2)}(z \p S_1 \m I_1)^{\beta+2} \bigr)\nonumber
\end{eqnarray}
for all $s \in[0,k]$ where the right-hand side defines an
integrable random variable since $\beta\p2 \in(0,\alpha)$. (Note
that without loss of generality, we can assume that $\beta$ is close
enough to $\alpha\m2$ so that $\beta\p2 > 0$.) Letting $n
\rightarrow\infty$ in (\ref{4.48}) and using the dominated
convergence theorem (twice) we see that (\ref{4.41}) holds as
claimed. This completes the proof.
\end{pf}

6. We now establish a remarkable probabilistic representation of
the global solution (\ref{4.31}) to the equation (\ref{4.25}). For
this, let us set
%
%
\begin{equation}
\label{4.50}
V_1(z) = \EE(z \vee S_1 \m X_1)^p
\end{equation}
for all $z \in\R_+$. From (\ref{3.13}), we see formally that $V_1(z)
= \EE_z e^{-p \sigma_\infty} G(Z_{\sigma_\infty})$ for all $z \in
\R_+$ where $\sigma_\infty= \inf\{ s \ge0 \dvtx Z_s = \infty
\}$, and this suggests that $z \mapsto V_1(z)$ should solve the
equation (\ref{4.25}) on $\R_+$. This can be derived rigourously as
follows.
\begin{Proposition}\label{prop8}
Let $F_1$ denote the global solution
(\ref{4.31}) to (\ref{4.25}) on~$\R_+$ with $F_1(0) =
1$. Then the following identity holds:
%
%
\begin{equation}
\label{4.51}
V_1(z) = a_1 F_1(z)
\end{equation}
for all $z \in\R_+$ where the constant $a_1$ is given
explicitly by
%
%
\begin{equation}
\label{4.52}
a_1 = \alpha(c \Gamma(-\alpha) )^{p/\alpha}
\frac{\Gamma(p)}{\Gamma(p/\alpha)} .
\end{equation}
\end{Proposition}
\begin{pf}
1. We first show that the identity (\ref{4.51})
holds with some cons\-tant $a_1\,{>}\,0$. For this, fix an arbitrary
$z_1\,{>},0$, set $F(z)\,{=}\,a F_1(z)$ for \mbox{$z\,{\in}\,\R_+$} where $a =
V_1(z_1)/F_1(z_1)$, and consider $\sigma_{ z_1} = \inf\{ s \ge
0 \dvtx Z_s = z_1 \}$. Then by (\ref{4.41}) and (\ref{4.25}), we find
that
%
%
\begin{equation}
\label{4.53}
F(z) = \EE_z e^{-p \sigma_{ z_1}} F(Z_{\sigma_{ z_1}}) =
F(z_1) \EE_z e^{-p \sigma_{ z_1}} = V_1(z_1) \EE_z
e^{-p \sigma_{ z_1}}
\end{equation}
for all $z \in[0,z_1]$. In addition, consider $\sigma_n = \inf
\{ s \ge0 \dvtx Z_s = n \}$ and set
%
%
\begin{equation}
\label{4.54}
V^n(z) = \EE_z e^{-p \sigma_n} G(Z_{\sigma_n})
\end{equation}
for $n > z_1$ and $z \in[0,z_1]$. Note that (\ref{3.13}) implies
that $V^n(z) \rightarrow V_1(z)$ as $n \rightarrow\infty$ for all
$z \in[0,z_1]$. Fixing $n > z_1$ and applying the strong Markov
property of $Z$ at $\sigma_{ z_1}$ we find that
%
%
\begin{eqnarray}
\label{4.55}
V^n(Z_{\sigma_{ z_1}}) &=& \EE_{Z_{\sigma_{ z_1}}} e^{-p
\sigma_n} G(Z_{\sigma_n}) \nonumber\\
&=& \EE_z \bigl( e^{-p \sigma_n
\circ\theta_{\sigma_{ z_1}} - p \sigma_{ z_1} + p
\sigma_{ z_1}} G(Z_{\sigma_n}) \circ\theta_{\sigma_{ z_1}}
\vert\cF_{t(\sigma_{ z_1})} \bigr) \\
&=& e^{p \sigma_{ z_1}} \EE_z \bigl( e^{-p \sigma_n}
G(Z_{\sigma_n}) \vert\cF_{t(\sigma_{ z_1})} \bigr)\nonumber
\end{eqnarray}
for all $z \in[0,z_1]$. Multiplying both sides by $e^{-p
\sigma_{ z_1}}$ and then taking $\EE_z$, we get
%
%
\begin{equation}
\label{4.56}
V^n(z_1) \EE_z e^{-p \sigma_{ z_1}} = V^n(z)
\end{equation}
for all $z \in[0,z_1]$ and $n > z_1$. Letting $n \rightarrow
\infty$ we obtain
%
%
\begin{equation}
\label{4.57}
V_1(z_1) \EE_z e^{-p \sigma_{ z_1}} = V_1(z)
\end{equation}
for all $z \in[0,z_1]$. Comparing (\ref{4.57}) with (\ref{4.53}), we
see that $V_1(z) = F(z)$ for all $z \in[0,z_1]$. Since $z_1>0$ was
arbitrary this establishes (\ref{4.51}) with some constant $a_1>0$.

2. To derive (\ref{4.52}), we may apply the Laplace transform $\LL$
on both sides of (\ref{4.25}) where $F(z) = V_1(z) = a_1 F_1 (z)$
for $z \in\R_+$ so that $a_1 = V_1(0)$. Using the well-known
properties (i)--(iii) recalled following (\ref{4.32}) above and
(iv)~$\LL[z F'(z)](\lambda) = -\lambda\LL[F]'(\lambda)
\m\LL[F](\lambda)$ for $\lambda>0$, it can be verified using
(\ref{4.28}) that this leads to
%
%
\begin{eqnarray}
\label{4.58}
&&\LL[F]'(\lambda) + \biggl( \frac{1 \p p}{\lambda} - \frac{c
\Gamma(2 \m\alpha)}{\alpha\m1} \lambda^{\alpha-1} \biggr)
\LL[F](\lambda) \nonumber\\[-8pt]\\[-8pt]
&&\qquad= - F(0) \frac{c
\Gamma(2 \m\alpha)}{\alpha\m1} \lambda^{\alpha-2}\nonumber
\end{eqnarray}
for $\lambda> 0$. Solving this equation under $\LL[F](\lambda)
\rightarrow0$ as $\lambda\rightarrow\infty$ [this condition is
satisfied since $F(z) = V_1(z) \sim z^p$ as $z \rightarrow\infty$
by (\ref{4.50}) above] we find that
%
%
\begin{equation}
\label{4.59}
\LL[F](\lambda) = \frac{F(0)}{(c \Gamma(-\alpha))^{p/\alpha}}
\frac{e^{c \Gamma(-\alpha) \lambda^\alpha}}{\lambda^{1
+p}} \Gamma\bigl(1 \p p/\alpha, c \Gamma(-\alpha)
\lambda^\alpha\bigr)
\end{equation}
for $\lambda>0$, where $\Gamma(a,x) = \int_x^\infty y^{a-1} e^{-y}
\,dy$ denotes the incomplete gamma function for $a>0$ and $x \ge0$.
Since $z \mapsto F(z)$ is increasing [by (\ref{4.50}) above], we can
use the Tauberian monotone density theorem (see, e.g., \cite{Ky},
Theorem~5.14, page 127) which states that (i) $\LL[F](\lambda) \sim
\ell\lambda^{-\rho}$ as $\lambda\downarrow0$ if and only if
(ii) $F(z) \sim(\ell/\Gamma(\rho)) z^{\rho-1}$ as $z \uparrow
\infty$ where $\rho>0$ and $\ell>0$. From (\ref{4.59}), we see that
(i) is satisfied with $\rho= 1 \p p$ and $\ell= (F(0)/(c
\Gamma(-\alpha))^{p/\alpha} ) \Gamma(1 \p p/\alpha)$ so that
(ii) yields (\ref{4.52}) since $F(z) = V_1(z) \sim z^p$ as $z
\rightarrow\infty$. This completes the proof.
\end{pf}

\section{Predicting the ultimate supremum}\label{sec5}

1. We will begin by connecting our findings on the free-boundary
problem from the previous section to the value function from
(\ref{3.2}).

\begin{Proposition}\label{prop9}
If the optimal stopping point $z_*$
from (\ref{3.14}) is finite, then the value function $V$ from
(\ref{3.2}) coincides on $[0,z_*]$ with $F$ from (\ref{4.31}) where
$a_0$ is set to $V(0)$. In terms of the function $V_1$ from
(\ref{4.50}), this reads as follows:
%
%
\begin{equation}
\label{5.1}
V(z) = a V_1(z)
\end{equation}
for all $z \in[0,z_*]$ where $a = V(0)/a_1 \in(0,1)$ and
$a_1$ is given by (\ref{4.52}) above. If the optimal stopping point
$z_*$ is not finite (i.e., the optimal stopping set~$D$ is empty),
then
%
%
\begin{equation}
\label{5.2}
V(z) = V_1(z)
\end{equation}
for all $z \in\R_+$.
\end{Proposition}
\begin{pf}
If $z_* < \infty$ then
%
%
\begin{equation}
\label{5.3}
V(z) = \EE_z e^{-p \sigma_{ z_*}} G(Z_{\sigma_{ z_*}})
= V(z_*) \EE_z e^{-p \sigma_{ z_*}}
\end{equation}
for all $z \in[0,z_*]$. Moreover, if we set $F(z) = a_0 F_1(z)$
for all $z \in\R_+$ with $a_0 = V(z_*)/F_1(z_*)$ then by
(\ref{4.41}) and (\ref{4.25}), we have
%
%
\begin{eqnarray}
\label{5.4}
F(z) &=& \EE_z e^{-p \sigma_{ z_*}} F(Z_{\sigma_{
z_*}})
= F(z_*) \EE_z e^{-p \sigma_{ z_*}} \nonumber\\[-8pt]\\[-8pt]
&=& V(z_*) \EE_z
e^{-p \sigma_{ z_*}}\nonumber
\end{eqnarray}
for all $z \in[0,z_*]$. Comparing (\ref{5.3}) and (\ref{5.4}), we
see that $V(z) = F(z)$ for all $z \in[0,z_*]$. Hence $a_0 = V(0)$
and this establishes (\ref{5.1}) upon recalling (\ref{4.51}). If
$z_*=\infty$ then (\ref{5.2}) follows from (\ref{3.13}) above. This
completes the proof.
\end{pf}

From (\ref{5.1}) and (\ref{5.2}), we see that the value function $V$
is a constant multiple of the function $V_1$ from (\ref{4.50}) up to
the first contact point with $G$ (when starting from 0 and moving
toward $\infty$ in the state space). The unknown constant needs to
be chosen so that the contact with $G$ occurs smoothly. Since $V \le
V_1$ this leads the following criterion for $D$ to be nonempty:
%
%
\begin{equation}
\label{5.5}
z_* < \infty\quad\mbox{if and only if}\quad \exists z_1 \in\R_+
\mbox{ such that } V_1(z_1) \ge G(z_1)
\end{equation}
or equivalently, the following criterion for $D$ to be empty:
%
%
\begin{equation}
\label{5.6}
z_* = \infty\quad\mbox{if and only if}\quad V_1(z) < G(z)
\mbox{ for all } z \in\R_+ .
\end{equation}
We will continue our analysis by examining when (\ref{5.5})
holds.\vadjust{\goodbreak}

2. Consider the function $H \dvtx [0,\infty) \rightarrow\R$ defined by
%
%
\begin{equation}
\label{5.7}
H(z) = (\LL_Z G \m p G)(z)
\end{equation}
for $z \ge0$ where $H(0) := H(0+)$ exists by (\ref{5.11}) below.
Recall that (\ref{4.41}) reads
%
%
\begin{equation}
\label{5.8}
\EE_z e^{-p \sigma} G(Z_\sigma) = G(z) + \EE_z \int_0^\sigma
e^{-ps} H(Z_s) \,ds
\end{equation}
for $z \in\R_+$ where $\sigma$ is any stopping time of $Z$ like in
Proposition \ref{prop7}. Set
%
%
\begin{equation}
\label{5.9}
N = \{ z \in[0,\infty) \dvtx H(z)<0 \} \quad\mbox{and}\quad P
= \{ z \in[0,\infty) \dvtx H(z) \ge0 \} .
\end{equation}
Then the following two inclusions are valid:
%
%
\begin{equation}
\label{5.10}
N \subseteq C \quad\mbox{and}\quad D \subseteq P .
\end{equation}

Indeed, to show the first inclusion (the second one then being
obvious) take any $z \in N$ and choose $\eps>0$ small enough such
that $(z \m\eps, z \p\eps) \cap\R_+ \subset N$ (note that $N$ is
open in $\R_+$). Inserting the stopping time $\sigma_\eps= \inf
\{ s \ge0 \dvtx Z_s \notin(z \m\eps, z \p\eps) \}$ into
(\ref{5.8}), we see that $\EE_z e^{-p \sigma_\eps}
G(Z_{\sigma_\eps})<$ $G(z)$ since $H(Z_s)<0$ for $s \in
[0,\sigma_\eps)$. Hence, $z$ belongs to $C$ as claimed.

3. Motivated by the important role that the function $H$ plays in
the optimal stopping problem (\ref{3.2}), we now determine its
asymptotic behavior at zero and infinity. Note that (\ref{5.11})
below and (\ref{5.10}) above imply (since~$H$ is continuous) that
the continuation set $C$ always contains the interval $[0,\eps)$ for
some $\eps>0$ sufficiently small so that the optimal stopping point
$z_*$ from~(\ref{3.14}) is always strictly larger than zero.

\begin{Proposition}\label{prop10}
The following relations are valid:
%
%
\begin{eqnarray}
\label{5.11}
\lim_{z \downarrow0} H(z) &=& -p G(0) = -p \EE S_1^p
< 0 ,\\
\label{5.12}
\lim_{z \uparrow\infty} z^{\alpha-p} H(z)
&=& \frac{c p}{\Gamma(p \m\alpha\p1)} \bigl( \Gamma(p \m
\alpha) \m\Gamma(p) \Gamma(1 \m\alpha) \bigr) .
\end{eqnarray}
\end{Proposition}
\begin{pf}
Since $G'(z) = p z^{p-1} F_{S_1}(z)$ and $G''(z)
= p (p \m1) z^{p-2} F_{S_1}(z) + p z^{p-1} f_{S_1}(z)$ we
see by (\ref{4.1}) and (\ref{4.5}) that
%
%
\begin{eqnarray}
\label{5.13}
H(z) &=& z G'(z) + \frac{c}{\alpha\m1} \int_0^z \frac{G''(x)}
{(z \m x)^{\alpha-1}} \,dx - p G(z) \nonumber\\
&=& p z^p
F_{S_1}(z) \nonumber\\[-8pt]\\[-8pt]
&&{}+ \frac{c}{\alpha\m1} \int_0^z \frac{p (p \m1)
x^{p-2} F_{S_1}(x) + p x^{p-1} f_{S_1}(x)}{(z \m x)^{\alpha-1}}
\,dx\hspace*{-28pt}\nonumber\\
&&{} - p G(z)\nonumber
\end{eqnarray}
for $z > 0$. Recall that the following asymptotic relations are
valid (see \cite{BDP}, Corollary 3):
%
%
\begin{eqnarray}
\label{5.14}
f_{S_1}(z) &\sim&\frac{z^{\alpha-2}}{(c \Gamma(-\alpha))^{1-1/
\alpha} \Gamma(\alpha\m1) \Gamma(1/\alpha)} \qquad\mbox{as }
z \downarrow0, \\
\label{5.15}
F_{S_1}(z) &\sim&\frac{z^{
\alpha-1}}{(c \Gamma(-\alpha))^{1-1/\alpha} \Gamma(\alpha)
\Gamma(1/\alpha)} \qquad\mbox{as } z \downarrow0 .
\end{eqnarray}
Using (\ref{5.14}) and (\ref{5.15}) together with (\ref{4.29}) it is
readily verified that the integral in (\ref{5.13}) tends to $0$ as
$z \downarrow0$. This easily yields the first equality in~%
(\ref{5.11}) and the second equality follows from (\ref{3.3}).

Moreover, using (\ref{2.12}) above we can further rewrite
(\ref{5.13}) as follows:
%
%
\begin{eqnarray}
\label{5.16}
H(z) &=& -p z^p \bigl(1 \m F_{S_1}(z) \bigr)\nonumber\\
&&{} + \frac{c p
(p \m1)} {\alpha\m1} \int_0^z \frac{x^{p-2}}{(z \m x)^{\alpha-1}}
\,dx\nonumber\\
&&{} - p \int_{z^p}^\infty\bigl(1 \m F_{S_1}(x^{1/p}) \bigr)
\,dx \\
&&{}
- \frac{c p (p \m1)} {\alpha\m
1} \int_0^z \frac{x^{p-2} (1 \m F_{S_1}(x))}{(z \m x)^{\alpha-1}}
\,dx \nonumber\\
&&{}+ \frac{c p}{\alpha\m1} \int_0^z \frac{x^{p-1} f_{S_1}
(x)}{(z \m x)^{\alpha-1}} \,dx\nonumber
\end{eqnarray}
for $z>0$. Recall that the following asymptotic relations are valid
(cf. \cite{BDP,Do,Pa}):
%
%
\begin{eqnarray}
\label{5.17}
f_{S_1}(z) &\sim&\frac{c}{z^{1+\alpha}} \qquad\mbox{as }
z \uparrow\infty,\\
\label{5.18}
1 \m F_{S_1}(z) &\sim&
\frac{c}{\alpha z^\alpha} \qquad\mbox{as } z \uparrow
\infty.
\end{eqnarray}
Using (\ref{5.17}) and (\ref{5.18}) together with (\ref{4.29}) it is
somewhat lengthy but still straightforward to verify that the final
two integrals in (\ref{5.16}) are $o(z^{p-\alpha})$ as $z
\rightarrow\infty$, whilst the first three terms in (\ref{5.16})
multiplied by $z^{\alpha-p}$ converge to the constant on the
right-hand side of (\ref{5.12}) as $z \rightarrow\infty$. This
completes the proof.
\end{pf}

4. Motivated by the identity (\ref{5.12}) let us consider the
function $\ell$ defined by
%
%
\begin{equation}
\label{5.19}
\ell(\alpha,p) = \frac{c p}{\Gamma(p \m\alpha\p1)}
\bigl( \Gamma(p \m\alpha) \m\Gamma(p) \Gamma(1 \m\alpha) \bigr)
\end{equation}
for $\alpha\in(1,2)$ and $p \in(1,\alpha)$. A direct examination
of the right-hand side in (\ref{5.19}) shows that there exist
$\alpha_* \in(1,2)$ (equal to $1.57$ approximately) and a strictly
increasing function $p_* \dvtx (\alpha_*,2) \rightarrow(1,2)$
satisfying\vadjust{\goodbreak} $p_*(\alpha_*+) = 1$, $p_*(2-) = 2$ and $p_*(\alpha) <
\alpha$ for $\alpha\in(\alpha_*,2)$ such that (i) $\ell(\alpha,p)
> 0$ if $\alpha\in(\alpha_*,2)$ and $p \in(1,p_*(\alpha))$; (ii)
$\ell(\alpha,p) < 0$ if either $\alpha\in(1,\alpha_*)$ and $p \in
(1,\alpha)$ or $\alpha\in[\alpha_*,2)$ and $p \in
(p_*(\alpha),\alpha)$; and (iii) $\ell(\alpha,p_*(\alpha)) = 0$ for
$\alpha\in(\alpha_*,2)$. Note that the properties (i)--(iii) do not
depend on the value of the constant~$c$ in~(\ref{2.4}). Recall also
from (\ref{5.12}) above that
%
%
\begin{equation}
\label{5.20}
\ell(\alpha,p) = \lim_{z \uparrow\infty} z^{\alpha-p} H(z)
\end{equation}
for all $\alpha\in(1,2)$ and $p \in(1,\alpha)$. In view of
(\ref{5.10}) this suggests that the sign of $\ell$ plays an
important role in the problem (\ref{3.2}).

Building on the facts presented in the previous sections, and
extending these arguments further in the proof below, we can now
present the main result of the paper. It should be recalled in the
statement below that the function $V_1$ can be expressed
probabilistically by (\ref{4.50}) and analytically by
(\ref{4.51})$+$(\ref{4.52}) [where $F_1$ is given by (\ref{4.31}) with
$a_0=1$], and the probabilistic and analytic representations of the
function $G$ are given in (\ref{2.12}) above (upon recalling that
$F_{S_1}$ admits an explicit series representation as shown in
\cite{BDP}, Theorem 1).
\begin{Theorem}\label{theo11}
\textup{I.} If $\alpha\in(\alpha_*,2)$ and $p
\in(1,p_*(\alpha))$ then there exists $z_* \in(0,\infty)$ such
that the stopping time $(\ref{3.14})$ is optimal in the problem
(\ref{3.2}) under $\PP_{ z}$ for $z \in[0,z_*]$. The optimal
stopping point $z_*$ can be characterized as the minimal $z \in
(0,\infty)$ for which
%
%
\begin{equation}
\label{5.21}
\beta_* V_1(z) \vert_{z=z_*} = G(z) \vert_{z=z_*},
\end{equation}
where $\beta_* \in(0,1)$ is the minimal $\beta\in(0,1)$ for
which (\ref{5.21}) has at least one root $z \in
(0,\infty)$. The optimal $z_*$ and $\beta_*$ satisfy the smooth fit
condition
%
%
\begin{equation}
\label{5.22}
\beta_* V_1'(z) \vert_{z=z_*} = G'(z) \vert_{z=z_*} .
\end{equation}
The value function from (\ref{3.2}) is given by $V(z) =
\beta_* V_1(z) = \beta_* \EE(z \vee S_1 \m X_1)^p$ for $z \in
[0,z_*]$.

\textup{II.} The stopping time (\ref{1.3}) is optimal in the problem
(\ref{2.3}) and the value from (\ref{2.3}) is given by $V =
T^{p/\alpha} \beta_* V_1(0) = T^{p/\alpha} \beta_* \EE(S_1 \m
X_1)^p = T^{p/\alpha} \beta_* \alpha(c\times\break
\Gamma(-\alpha))^{p/\alpha} \Gamma(p)/\Gamma(p/\alpha)$.
\end{Theorem}
\begin{pf}
Since part II follows from part I as discussed in
Sections \ref{sec2} and~\ref{sec3} above, it is enough to prove part I. For this, we
will first show that the assumptions $\alpha\in(\alpha_*,2)$ and
$p \in(1,p_*(\alpha))$ imply the existence of $z_1>0$ (large
enough) such that
%
%
\begin{equation}
\label{5.23}
V_1(z) > G(z)
\end{equation}
for all $z \ge z_1$. We will then show how the knowledge of
(\ref{5.23}) combined with the properties and facts about $V_1$ and
$G$ derived in the previous sections yield the existence of
$\beta_*$ and $z_*$ satisfying the remaining statements\vadjust{\goodbreak} of part~I.

1. To prove (\ref{5.23}), recall that the identity (\ref{4.41}) is
applicable to $G$ in place of $F$ with $\sigma\equiv n$ for $n \ge
1$. Letting $n \rightarrow\infty$ in this identity, using
(\ref{3.13}) combined with the fact that each $e^{-pn} G(Z_n^z)$ is
dominated by $(z \vee S_1 \m I_1)^p + \EE S_1^p$ which clearly
has finite expectation, as well as the fact that the function~$H$ is
bounded (by the result of Proposition \ref{prop10}), it follows by the
dominated convergence theorem that
%
%
\begin{equation}
\label{5.24}
\EE(z \vee S_1 \m X_1)^p = G(z) + \EE\int_0^\infty e^{-ps}
H(Z_s^z) \,ds
\end{equation}
for all $z \ge0$. Recognizing the left-hand side of (\ref{5.24}) as
$V_1(z)$, we see that~(\ref{5.23}) will be established if we show the
existence of $z_1>0$ (large enough) such that
%
%
\begin{equation}
\label{5.25}
I(z) := \EE\int_0^\infty e^{-ps} H(Z_s^z) \,ds > 0
\end{equation}
for all $z \ge z_1$.

To show (\ref{5.25}) recall from (i) following (\ref{5.19}) above
that $\ell:= \ell(\alpha,p)$ in~(\ref{5.20}) is strictly positive
when $\alpha\in(\alpha_*,2)$ and $p \in(1,p_*(\alpha))$ are given
and fixed. Hence for any given and fixed $\eps>0$ (small) there
exists $z_\eps>0$ (large) such that
%
%
\begin{equation}
\label{5.26}
z^{\alpha-p} H(z) \ge\ell\m\eps
\end{equation}
for all $z \ge z_\eps$. Consider
%
%
\begin{eqnarray}
\label{5.27}
J(z) &:=& \EE\int_0^\infty e^{-ps} H(Z_s^z) I(Z_s^z < z_\eps)
\,ds, \\
K(z) &:=& \EE\int_0^\infty e^{-ps} H(Z_s^z) I(Z_s^z \ge
z_\eps) \,ds
\end{eqnarray}
and note that $I(z) = J(z) \p K(z)$ for all $z \ge0$.

Let $M>0$ be large enough so that $\vert H(z) \vert\le M$ for all
$z \ge0$. Then we have
%
%
\begin{equation}
\label{5.29}
\vert J(z) \vert\le M \int_0^\infty e^{-ps} \PP(Z_s^z < z_\eps) \,ds
\end{equation}
for all $z \ge0$. Moreover, by (\ref{5.18}) we see that
%
%
\begin{eqnarray}
\label{5.30}
\PP(Z_s^z < z_\eps) &=& \PP\bigl( e^s\bigl(z \vee S_{t(s)} \m X_{t(s)}\bigr)
< z_\eps\bigr) \nonumber\\
&\le&\PP\bigl( z \vee S_{t(s)} \m S_{t(s)} < z_\eps
\bigr) \nonumber\\
&\le&\PP\bigl( z \m S_{t(s)} < z_\eps\bigr)\\
&\le&
\PP( S_1 > z \m z_\eps) \nonumber\\
&\le& N \frac{c}{\alpha} (z \m
z_\eps)^{-\alpha}\nonumber
\end{eqnarray}
for all $z > z_\eps$ with some $N>0$ large enough. Combining
(\ref{5.29}) and (\ref{5.30}) we find that
%
%
\begin{equation}
\label{5.31}
\vert J(z) \vert\le\frac{M N c}{p \alpha} (z \m
z_\eps)^{-\alpha}
\end{equation}
for all $z > z_\eps$.

On the other hand, by (\ref{5.26}) we see that
%
%
\begin{eqnarray}
\label{5.32}\qquad
K(z) &=& \EE\int_0^\infty e^{-ps} H(Z_s^z) I(Z_s^z
\ge z_\eps) \,ds \nonumber\\
&\ge&(\ell\m\eps) \int_0^\infty e^{-ps} \EE[
(Z_s^z)^{p-\alpha}I(Z_s^z \ge z_\eps)] \,ds \nonumber\\
&=& (\ell
\m\eps) \int_0^\infty e^{-\alpha s} \EE\bigl[ \bigl(z \vee S_{t(s)}
\m X_{t(s)}\bigr)^{p-\alpha}I(Z_s^z \ge z_\eps)\bigr] \,ds
\nonumber\\[-8pt]\\[-8pt]
&=& (\ell\m\eps) z^{p-\alpha} \int_0^\infty e^{-\alpha s} \EE\biggl[
\biggl(1 \vee\frac{S_{t(s)}}{z} \m\frac{X_{t(s)}}{z} \biggr)^{p
-\alpha}I(Z_s^z \ge z_\eps)\biggr] \,ds \nonumber\\
&\ge&(\ell\m\eps)
z^{p-\alpha} \int_0^\infty e^{-\alpha s} \EE[(1 \vee
S_1 \m I_1)^{p-\alpha}I(Z_s^z \ge z_\eps)] \,ds \nonumber\\ &
\ge&\frac{(\ell\m\eps)}{\alpha} \bigl( \EE(1 \vee S_1 \m I_1)^{
p-\alpha} - \delta\bigr) z^{p-\alpha}\nonumber
\end{eqnarray}
for all $z \ge1 \vee z_\delta$, where in the second last inequality
we use that
%
%
\begin{equation}
\label{5.33}
1 \vee\frac{S_{t(s)}}{z} - \frac{X_{t(s)}}{z} \le1 \vee
\frac{S_1}{z} - \frac{I_1}{z} \le1 \vee S_1 - I_1
\end{equation}
for all $s \ge0$ and $z \ge1$, and in the last inequality we use
that
%
%
\begin{eqnarray}
\label{5.34}
&&\lim_{z \rightarrow\infty} \int_0^\infty e^{-\alpha s} \EE[
(1 \vee S_1 \m I_1)^{p-\alpha}I(Z_s^z \ge z_\eps)] \,ds
\nonumber\\[-8pt]\\[-8pt]
&&\qquad=
\frac{1}{\alpha} \EE(1 \vee S_1 \m I_1)^{p-\alpha} < \infty\nonumber
\end{eqnarray}
by the dominated convergence theorem since $Z_s^z \rightarrow
\infty$ as $z \rightarrow\infty$ [from (\ref{5.34}) we see that for
given $\delta\in(0,\EE(1 \vee S_1 \m I_1)^{p-\alpha})$ there
exists $z_\delta>0$ such that the final inequality in (\ref{5.32})
holds for all $z \ge z_\delta$]. Since the right-hand side in
(\ref{5.31}) tends faster to zero than the right-hand side in
(\ref{5.32}) as $z \uparrow\infty$, we see that (\ref{5.23}) holds
with some $z_1>0$ large enough as claimed.

2. We now establish the existence of $\beta_*$ and $z_*$ satisfying
the remaining statements of part I. For this, recall that
(\ref{5.23}) holds for $z=z_1$ so that for some $\beta_1 \in(0,1)$
sufficiently close to $1$ we have $\beta_1 V_1(z_1) > G(z_1)$. Since
$\beta_1 V_1(z) \sim\beta_1 z^p < z^p \sim G(z)$ as $z \rightarrow
\infty$ we also see that there exists $z_2 > z_1$ such that $\beta_1
V_1(z) < G(z)$ for all $z \ge z_2$. This shows that for some
$\beta_0 \in(0,1)$ sufficiently close to $0$ we have $\beta_0
V_1(z) < G(z)$ for all $z \ge0$ [recall that $V_1(0) = \EE(S_1 -
X_1)^p > 0$ and that $V_1$ is increasing]. It follows therefore by
continuity that there exists the smallest $\beta_* \in
(\beta_0,\beta_1) \subset(0,1)$ such that the set $A = \{ z \in
\R_+ \vert\beta_* V_1(z) = G(z) \}$ is nonempty so that
$\beta V_1(z) < G(z)$ for all $z \in\R_+$ if $\beta\in
(0,\beta_*)$. Setting $w_* = \inf A$ we see that $w_*$ belongs to
$A$ by continuity so that (\ref{5.21}) holds for $z=w_*$. Moreover,
since $V_1(z_1) > G(z_1)$ we know by (\ref{5.5}) that $z_* = \inf D
< \infty$ so that by (\ref{5.1}) we have $V(z) = a_* V_1(z)$ for all
$z \in[0,z_*]$ with some $a_* \in(0,1)$. By the construction of
$\beta_*$ and $w_*$ it follows therefore that $\beta_* \le a_*$ and
$w_* \ge z_*$. If either $\beta_*<a_*$ or equivalently $w_*>z_*$,
then since $\beta_* V_1(z) = \EE e^{-p \sigma_{w_*}}(\beta_*
V_1)(Z_{\sigma_{w_*}}^z)$ for all $z \in[0,w_*]$ by the result of
Proposition \ref{prop7}, and this further equals $\EE e^{-p
\sigma_{w_*}}G(Z_{\sigma_{w_*}}^z)$ for all $z \in[0,w_*]$ by
definition of $\sigma_{w_*}$, we see that $\beta_* V_1(0) \ge V(0)$
while at the same time $\beta_* V_1(0) < a_* V_1(0) = V(0)$ which is
a contradiction. Thus $\beta_* = a_*$ and $w_* = z_*$ so that $V(z)
= \beta_* V_1(z)$ for all $z \in[0,z_*]$ as claimed. The smooth fit
condition (\ref{5.22}) then follows by the result of Proposition \ref{prop5}.
This completes the proof of part~I whence part II follows as
discussed above.
\end{pf}

5. In the final part of this section, we briefly consider the case
when the hypotheses of Theorem \ref{theo11} are not satisfied.
\begin{Proposition}\label{prop12}
If either $\alpha\in(1,\alpha_*)$
or $p
\in(p_*(\alpha),\alpha)$, then there exists $z_1 > 0$ large enough
such that $V_1(z) < G(z)$ for all $z \ge z_1$.
\end{Proposition}
\begin{pf}
This can be proved in exactly the same way as
(\ref{5.23}) above upon noting that $\ell:= \ell(\alpha,p)$ in
(\ref{5.20}) is strictly negative when either $\alpha\in
(1,\alpha_*)$ or $p \in(p_*(\alpha),\alpha)$ and replacing
(\ref{5.26}) with $z^{\alpha-p} H(z) \le\ell\p\eps$ for all $z
\ge z_\eps$. This leads to (\ref{5.31}) without changes and
(\ref{5.32}) holds with the inequalities reversed since $\ell\p
\eps< 0$ in this case. Different rates of convergence in the
resulting inequalities then complete the proof just as above.
\end{pf}

It follows from the result of Proposition \ref{prop12} that the continuation
set~$C$ contains the interval $[z_1,\infty)$ for some $z_1>0$ large
enough when either $\alpha\in(1,\alpha_*)$ or $p \in(p_*(\alpha),
\alpha)$. It shows that the stopping time (\ref{1.3}) can no longer
be optimal in this case (in the sense that it is not optimal to stop
at $t \in[0,T)$ when $S_t \m X_t$ is sufficiently large). This
stands in sharp contrast with the Brownian motion case (formally
corresponding to $\alpha=2$) where it is optimal to stop in such a
case. Recall also that the continuation set~$C$ always contains the
interval $[0,\eps)$ for some $\eps>0$ sufficiently small so that the
stopping set~$D$ must be contained in $[\eps,z_1 \m\delta]$ for
some $\delta>0$. We do not know whether $V_1(z) < G(z)$ holds for
all $z \in\R_+$ in this case, or equivalently, whether the stopping
set $D$ is empty [recall (\ref{5.6}) above]. This is an interesting
open question. We refer to \cite{DP-1}, Figure 1, for a related
phenomenon in the presence of strictly positive drifts and the
absence of jumps.

\begin{appendix}\label{app}
\section*{Appendix}

In this section, we determine the action of the infinitesimal
generator of the reflected process $Y=S-X$ when $X$ is a\vadjust{\goodbreak} general
(strictly) stable L\'evy process (see \cite{Wa}). Set
%
%
\setcounter{equation}{0}
\begin{equation}
\label{A.1}
\nu_\alpha(dx) = \frac{c_+}{x^{1+\alpha}} I(x > 0) \,dx +
\frac{c_-}{(-x)^{1+\alpha}} I(x < 0) \,dx,
\end{equation}
where $c_+$ and $c_-$ are nonnegative constants (not both zero) and
$\alpha\in(0,2)$. For $\alpha=1$, the two constants need to be
identical (see, e.g., \cite{Sa}, pages 86 and 87), so that
%
%
\begin{equation}
\label{A.2}
\nu_1(dx) = \frac{c}{x^2} I(x \ne0) \,dx
\end{equation}
with $c>0$. Recall that $C_b^2(\R_+)$ denotes the class of twice
continuously differentiable functions $F \dvtx \R_+ \rightarrow\R$ such
that $F'$ and $F''$ are bounded on $\R_+$.

\begin{Proposition}\label{propA1}
Let $X=(X_t)_{t \geq0}$ be a stable
L\'evy process of index $\alpha\in(1,2)$ whose characteristic
function is given by
%
%
\begin{eqnarray}
\label{A.3}
\EE e^{i \lambda X_t} &=& \exp\biggl(t \int_{-\infty}^\infty(e^{i
\lambda x} - 1 - i \lambda x) \nu_\alpha(dx) \biggr) \nonumber\\[-8pt]\\[-8pt]
&=& e^{(c_+(-i
\lambda)^\alpha+ c_- (i \lambda)^\alpha) \Gamma(-\alpha) t}\nonumber
\end{eqnarray}
for $\lambda\in\R$ and $t \ge0$. Then the infinitesimal
generator $\LL_Y$ of the reflected process $Y = S \m X$ takes any of
the following three forms for $y>0$ given and fixed:
%
%
\begin{itemize}
\item[]\mbox{It\^o's form}
\begin{eqnarray}
\label{A.4}
\LL_Y F(y) &=& \int_0^y
\bigl( F(y \m x) \m F(y) \p F'(y) x \bigr) \frac{c_+}{x^{1+\alpha}}
\,dx \nonumber\\
&&{}+ \frac{c_+ (F(0) \m F(y))}{\alpha y^\alpha} + \frac{c_+
F'(y)}{(\alpha\m1) y^{\alpha-1}}\\
&&{}
+
\int_0^\infty\bigl( F(y \p x) \m F(y) \m F'(y) x \bigr)
\frac{c_-}{x^{1+\alpha}} \,dx ,\nonumber
\end{eqnarray}
\item[]\mbox{Riemann--Liouville's form}
\begin{eqnarray}
\label{A.5}
\LL_Y F(y) &=& \frac{c_+}
{\alpha(\alpha\m1)} \frac{d^2}{dy^2} \int_0^y \frac{F(x)}
{(y \m x)^{\alpha-1}} \,dx + \frac{c_+ F(0)}{\alpha y^\alpha}
\nonumber
\\[-8pt]
\\[-8pt]
\nonumber
&&{}
+ \frac{c_-}{\alpha(\alpha\m1)} \frac{d^2}{dy^2} \int_y^\infty
\frac{F(x)}{(x \m y)^{\alpha-1}} \,dx ,
\end{eqnarray}
\item[]\mbox{Caputo's form}
\begin{eqnarray}
\label{A.6}
\LL_Y F(y) &=& \frac{c_+}{\alpha
(\alpha\m1)} \int_0^y \frac{F''(x)}{(y \m x)^{\alpha- 1}}
\,dx
\nonumber
\\[-8pt]
\\[-8pt]
\nonumber
&&{}+ \frac{c_-}{\alpha(\alpha\m1)} \int_y^\infty
\frac{F''(x)}{(x \m y)^{\alpha- 1}} \,dx,\nonumber
\end{eqnarray}
\end{itemize}
whenever $F \in C_b^2(\R_+)$ satisfies
%
%
\begin{equation}
\label{A.7}
F'(0+) = 0 \qquad\mbox{(normal reflection)}
\end{equation}
with $\vert F''(y) \vert= O(y^\gamma)$ as $y \rightarrow
\infty$ for some $\gamma< \alpha\m2$ [as well as $\vert F(y)
\vert= O(y^\delta)$ as $y \rightarrow\infty$ for some $\delta<
\alpha\m2$ in (\ref{A.5}) above].
\end{Proposition}
\begin{pf}
As in the proof of Proposition \ref{prop4}, it is enough to
derive (\ref{A.4}). For this, fix $t>0$ and note that by It\^o's
formula we have
%
%
\begin{eqnarray}
\label{A.8}
F(Y_t) &=& F(Y_0) + \int_0^t F'(Y_{s-}) (dS_s - dX_s) \nonumber\\[-8pt]\\[-8pt]
&&{}+ \sum_{0
< s \le t} \bigl(F(Y_s) - F(Y_{s-}) - F'(Y_{s-}) (\Delta S_s
- \Delta X_s) \bigr)\nonumber
\end{eqnarray}
since $[Y,Y]^c \equiv0$ for the same reasons as in (\ref{4.7}).
Letting $S_s = S^c_s + S^d_s$ be the decomposition of $s \mapsto
S_s$ into continuous and discontinuous parts, and noting that
$dS^d_s = \Delta S_s$, we see that (\ref{A.8}) simplifies to
%
%
\begin{eqnarray}
\label{A.9}
F(Y_t) &=& F(Y_0) + M_t + \int_0^t F'(Y_{s-}) \,dS^c_s \nonumber\\[-8pt]\\[-8pt]
&&{}+ \sum_{0
< s \le t} \bigl(F(Y_{s-} \p\Delta Y_s) - F(Y_{s-}) + F'(Y_{s-})
\Delta X_s \bigr),\nonumber
\end{eqnarray}
where $M_t = - \int_0^t F'(Y_{s-}) \,dX_s$. Since $F'$ is bounded
the same argument as following (\ref{4.8}) above shows that $M$ is a
martingale. If $s$ belongs to the support of $dS^c_s$ in $[0,t]$,
then either $S_{s-\eps}^c < S_s^c$ and therefore $S_{s-\eps} < S_s$
for $\eps>0$ implying $Y_{s-}=0$, or $S_s^c < S_{s+\eps}^c$ and
therefore $S_s < S_{s+\eps}$ for $\eps>0$ implying $Y_s=0$. Since
there could be at most countably many $s$ in $[0,t]$ for which $Y_s
\ne Y_{s-}$, it follows using (\ref{A.7}) that the integral with
respect to~$dS^c_s$ in~(\ref{A.9}) is zero. Moreover, the right-hand
side of (\ref{A.9}) can further be rewritten as follows:
%
%
\begin{eqnarray}
\label{A.10}\qquad
F(Y_t) &=& F(Y_0) + M_t \nonumber\\
&&{}+ \sum_{0 < s \le t} \bigl( [ F(Y_{s-}
- \Delta X_s) \m F(Y_{s-}) \nonumber\\[-8pt]\\[-8pt]
&&\hspace*{88.4pt}{}\p F'(Y_{s-}) \Delta X_s ]
I(\Delta X_s \le Y_{s-}) \nonumber\\
&&\hspace*{40.1pt}{}
+ [F(0)
\m F(Y_{s-}) + F'(Y_{s-}) \Delta X_s ] I(\Delta X_s
> Y_{s-}) \bigr)\nonumber
\end{eqnarray}
using the same arguments as in (\ref{4.9}) above. Taking $\EE_y$ on
both sides of~(\ref{A.10}), where $\PP_{ y}$ denotes a probability
measure under which $Y_0=y$, and applying\vadjust{\goodbreak} the compensation formula
(see, e.g., \cite{RY}, page 475) we find that\vspace*{-2pt}
%
%
\begin{eqnarray}
\label{A.11}
&&\EE_y F(Y_t) \m F(y) \nonumber\\[-3pt]
&&\qquad= \EE_y \biggl[\int_0^t ds \biggl(\int_{-
\infty}^{Y_s} [ F(Y_s \m x) \m F(Y_s) \p F'(Y_s) x ]
\nu_\alpha(dx) \\[-3pt]
&&\qquad\quad\hspace*{62.3pt}{}
+ \int_{Y_s}^\infty
[ F(0) \m F(Y_s) \p F'(Y_s) x ] \nu_\alpha(dx)
\biggr) \biggr]\nonumber\vspace*{-2pt}
\end{eqnarray}
for all $y>0$. The applicability of this formula (see, e.g.,
\cite{Ky}, page 97) follows from the facts that $\vert F'(y) \vert
\le C$ and $\vert F''(y) \vert\le C$ for all $y \ge0$ with some
$C>0$ so that the mean value theorem yields the existence of
$\xi_{s,x} \in(Y_s,Y_s \p x)$ and $\eta_{s,x} \in(Y_s,Y_s \p x)$
such that\vspace*{-2pt}
%
%
\begin{eqnarray}
\label{A.12}\qquad
&&\EE_y \biggl[\int_0^t ds \biggl(\int_{-\infty}^0 \vert F(Y_s \m
x) \m F(Y_s) \p F'(Y_s) x \vert\nu_\alpha(dx) \biggr)
\biggr] \nonumber\\[-3pt]
&&\qquad
= \EE_y \biggl[\int_0^t ds \biggl(
\int_0^\infty\vert F(Y_s \p x) \m F(Y_s) \m F'(Y_s) x
\vert\frac{c}{x^{1+\alpha}} \,dx \biggr) \biggr] \nonumber\\[-3pt]
&&\qquad
\le\EE_y \biggl[ \int_0^t ds \biggl( \int_0^1 \frac{1}{2}
\vert F''(\xi_{s,x}) \vert x^2 \frac{c}{x^{1+\alpha}}
\,dx\\[-2pt]3
&&\hspace*{81.8pt}{} + \int_1^\infty\bigl( \vert F'(\eta_{s,x}) \vert x + \vert
F'(Y_s) \vert x \bigr) \frac{c}{x^{1+\alpha}} \,dx\biggr)\biggr] \nonumber\\[-3pt]
&&\qquad\le c \EE_y \biggl[ \int_0^t ds \biggl(\frac{C}{2(2-
\alpha)} + \frac{2C}{\alpha-1} \biggr) \biggr] = \frac{c (7-
3\alpha) C}{2(2-\alpha)(\alpha-1)} t < \infty,\nonumber\vspace*{-2pt}
\end{eqnarray}
where the remaining two integrals (from $0$ to $Y_s$ and from $Y_s$
to $\infty$) can be controlled (bound from above) in exactly the
same way as in (\ref{4.11}) above. Dividing both sides of
(\ref{A.11}) by $t$, letting $t \downarrow0$ and using the
dominated convergence theorem, we get\vspace*{-2pt}
%
%
\begin{eqnarray}
\label{A.13}\quad
\LL_Y F(y) &=& \int_{-\infty}^y [ F(y \m x) \m F(y) \p F'(y)
x ] \nu_\alpha(dx) \nonumber\\[-10pt]\\[-10pt]
&&{}
+ [ F(0) \m F(y)
] \int_y^\infty\nu_\alpha(dx) + F'(y) \int_y^\infty x
\nu_\alpha(dx),\nonumber\vspace*{-2pt}
\end{eqnarray}
which is easily verified to be equal to the right-hand side of
(\ref{A.4}) for all $y>0$ upon using (\ref{A.1}). This completes the
proof.\vspace*{-3pt}
\end{pf}

\begin{Proposition}\label{propA2}
Let $X=(X_t)_{t \geq0}$ be a stable
L\'evy process of index $\alpha\in(0,1)$ whose characteristic
function is given by
%
%
\begin{eqnarray}
\label{A.14}
\EE e^{i \lambda X_t} &=& \exp\biggl(t \int_{-\infty}^\infty(e^{i
\lambda x} - 1 ) \nu_\alpha(dx) \biggr)\nonumber\\[-10pt]\\[-10pt]
&=& e^{(c_+(-i \lambda)^\alpha
+ c_- (i \lambda)^\alpha) \Gamma(-\alpha) t}\nonumber\vadjust{\goodbreak}
\end{eqnarray}
for $\lambda\in\R$ and $t \ge0$. Then the infinitesimal
generator $\LL_Y$ of the reflected process $Y = S \m X$ takes any of
the following three forms for $y>0$ given and fixed:
%
%
\begin{itemize}
\item[]\mbox{It\^o's form}
\begin{eqnarray}
\label{A.15}
\LL_Y F(y)
&=& \int_0^y \bigl(F(y \m x) \m F(y) \bigr) \frac{c_+}{x^{1+\alpha}} \,dx
+
\frac{c_+ (F(0) \m F(y))}{\alpha y^\alpha}
\nonumber
\\[-8pt]
\\[-8pt]
\nonumber
&&{}+ \int_0^\infty\bigl( F(y \p x) \m F(y)\bigr) \frac{c_-}
{x^{1+\alpha}} \,dx,
\end{eqnarray}
\item[]\mbox{Riemann--Liouville's form}
\begin{eqnarray}
\label{A.16}
\LL_Y F(y)& =&
-\frac{c_+}{\alpha} \,\frac{d}{dy} \int_0^y \frac{F(x)}{(y \m x)^{
\alpha}} \,dx + \frac{c_+ F(0)}{\alpha y^\alpha}
\nonumber
\\[-8pt]
\\[-8pt]
\nonumber
&&{}+ \frac{
c_-}{\alpha} \,\frac{d}{dy} \int_y^\infty\frac{F(x)}{(x \m y)^{
\alpha}} \,dx,
\end{eqnarray}
\item[]\mbox{Caputo's form}
\begin{equation}
\label{A.17}
\LL_Y F(y) = -\frac{c_+}{\alpha} \int_0^y
\frac{F'(x)}{(y \m x)^{\alpha}} \,dx
+ \frac{c_-}{\alpha}
\int_y^\infty\frac{F'(x)}{(x \m y)^{\alpha}} \,dx,
\end{equation}
\end{itemize}
whenever $F \in C_b^2(\R_+)$ satisfies $\vert F'(y) \vert=
O(y^\gamma)$ as $y \rightarrow\infty$ for some $\gamma< \alpha\m
1$ [as well as $\vert F(y) \vert= O(y^\delta)$ as $y \rightarrow
\infty$ for some $\delta< \alpha\m1$ in (\ref{A.16}) above].
\end{Proposition}
\begin{pf}
As in the proof of Proposition \ref{prop4} it is enough to
derive (\ref{A.15}). For this, fix $t>0$ and note that since $X$ is
a pure jump semimartingale with bounded variation, we have $dX_s =
\Delta X_s$ and $dS_s = \Delta S_s$ for $0 < s \le t$, so that
It\^o's formula yields
%
%
\begin{equation}
\label{A.18}
F(Y_t) = F(Y_0) + \sum_{0 < s \le t} \bigl(F(Y_s) \m F(Y_{s-})
\bigr) .
\end{equation}
Proceeding as in (\ref{A.10}), taking $\EE_y$ on both sides of the
resulting identity and applying the compensation formula (see, e.g.,
\cite{RY}, page 475), we find that
%
%
\begin{eqnarray}
\label{A.19}
&&
\EE_y F(Y_t) \m F(y)\nonumber\\
&&\qquad= \EE_y \biggl[ \int_0^t ds \biggl(
\int_{-\infty}^{Y_s} [ F(Y_s \m x) \m F(Y_s) ]
\nu_\alpha(dx) \\
&&\qquad\quad\hspace*{63.4pt}{}
+ \int_{Y_s}^\infty
[ F(0) \m F(Y_s)] \nu_\alpha(dx) \biggr) \biggr]\nonumber
\end{eqnarray}
for all $y>0$. The applicability of this formula (see, e.g.,
\cite{Ky}, page 97) follows from the facts that $\vert F(y) \vert\le
C$ and $\vert F'(y) \vert\le C$ for all $y \ge0$ with some $C>0$
so\vadjust{\goodbreak} that the mean value theorem yields the existence of $\xi_{s,x}^1
\in(Y_s,Y_s \p x)$, $\xi_{s,x}^2 \in(Y_s \m x,Y_s)$ and $\eta_s
\in(0,Y_s)$ such that
%
%
\begin{eqnarray}
\label{A.20}
&&\EE_y \biggl[\int_0^t ds \biggl(\int_{-\infty}^{Y_s} \vert
F(Y_s \m x) \m F(Y_s) \vert\nu_\alpha(dx)\nonumber\\[-2pt]
&&\qquad\quad\hspace*{29.5pt}{} + \int_{Y_s}^\infty
\vert F(0) \m F(Y_s) \vert\nu_\alpha(dx) \biggr)
\biggr] \nonumber\\[-2pt]
&&\qquad
\le\EE_y \biggl[\int_0^t ds
\biggl(\int_0^\infty\vert F(Y_s \p x) \m F(Y_s) \vert
\frac{c}{x^{1+\alpha}} \,dx \nonumber\\[-2pt]
&&\hspace*{47.6pt}\qquad\quad{}
+ \int_0^{Y_s}
\vert F(Y_s \m x) \m F(Y_s) \vert\frac{c}{x^{1+\alpha}}
\,dx
\nonumber\\[-2pt]
&&\qquad\quad\hspace*{71.2pt}{}
+ \int_{Y_s}^\infty\vert F(0) \m F(Y_s) \vert
\frac{c}{x^{1+\alpha}} \,dx \biggr) \biggr] \nonumber\\[-2pt]
&&\qquad
\le\EE_y \biggl[\int_0^t ds \biggl( \int_0^1 \vert F'(\xi_{s,x}^1)
\vert x \frac{c}{x^{1+\alpha}} \,dx\nonumber\\[-2pt]
&&\qquad\quad\hspace*{48pt}{} + \int_1^\infty\vert
F(Y_s \p x) \m F(Y_s) \vert\frac{c}{x^{1+\alpha}} \,dx
\\[-2pt]
&&\qquad\quad\hspace*{48pt}{}
+ \int_0^1 \vert F'(\xi_{s,x}^2) \vert
x \frac{c}{x^{1+\alpha}} \,dx\nonumber\\[-2pt]
&&\qquad\quad\hspace*{48pt}{} + \int_1^\infty\vert F(Y_s
\m x) \m F(Y_s) \vert\frac{c}{x^{1+\alpha}} \,dx
\nonumber\\[-2pt]
&&\qquad\quad\hspace*{48pt}{}
+ \int_{Y_s}^1 \vert F'(\eta_s) \vert Y_s \frac{c}
{x^{1+\alpha}} \,dx\, I(Y_s \le1)\nonumber\\[-2pt]
&&\qquad\quad\hspace*{72pt}{} + \int_1^\infty
\vert F(0) \m F(Y_s) \vert\frac{c}{x^{1+\alpha}} \,dx
\biggr) \biggr] \nonumber\\[-2pt]
&&\qquad
\le c \EE_y \biggl[\int_0^t
ds \biggl( \frac{2C}{1-\alpha} + \frac{6C}{\alpha} + \frac{C}
{\alpha} (Y_s^{1-\alpha} \m Y_s ) I(Y_s \le1)
\biggr) \biggr] \nonumber\\[-2pt]
&&\qquad\le c \biggl( \frac{2C}{1-\alpha} + \frac{7C}
{\alpha} \biggr) t < \infty\nonumber
\end{eqnarray}
upon using that $1 \m\alpha\in(0,1)$ in the final inequality.
Dividing both sides of~(\ref{A.19}) by $t$, letting $t \downarrow0$
and using the dominated convergence theorem, we get
%
%
\begin{eqnarray}
\label{A.21}
\LL_Y F(y) &=& \int_{-\infty}^y [ F(y \m x) \m F(y) ]
\nu_\alpha(dx) \nonumber\\[-10pt]\\[-10pt]
&&{} + [ F(0) \m F(y) ] \int_y^\infty
\nu_\alpha(dx),\nonumber
\end{eqnarray}
which is easily verified to be equal to the right-hand side of
(\ref{A.15}) for all $y>0$ upon using (\ref{A.1}). This completes
the proof.
\end{pf}

\begin{Proposition}\label{propA3}
Let $X=(X_t)_{t \geq0}$ be a stable
L\'evy process of index $1$ whose characteristic function is given
by
%
%
\begin{equation}
\label{A.22}\qquad
\EE e^{i \lambda X_t} = \exp\biggl( t \int_{-\infty}^\infty\bigl(
e^{i \lambda x} - 1 - i \lambda x I(\vert x \vert\le
1) \bigr) \nu_1(dx) \biggr) = e^{-c \vert\lambda\vert
\pi t}
\end{equation}
for $\lambda\in\R$ and $t \ge0$. Then the infinitesimal
generator $\LL_Y$ of the reflected process $Y = S - X$ takes any of
the following three forms for $y>0$ given and fixed:
%
%
\begin{itemize}
\item[]\mbox{It\^o's form}\vspace*{-1pt}
\begin{eqnarray}
\label{A.23}
\LL_Y F(y) &=&
\int_0^y \bigl( F(y \m x) \m F(y) \p F'(y) x \bigr) \frac{c}{x^2}
\,dx \nonumber\\
&&{}+ \frac{c (F(0) \m F(y))}{y}
\nonumber
\\[-8pt]
\\[-8pt]
\nonumber
&&{}
+ \int_0^y
\bigl( F(y \p x) \m F(y) \m F'(y) x \bigr) \frac{c}{x^2} \,dx \\
&&{}+
\int_y^\infty\bigl( F(y \p x) \m F(y) \bigr) \frac{c}{x^2}
\,dx,\nonumber
\end{eqnarray}
\item[]\mbox{Riemann--Liouville's form}
\begin{equation}
\label{A.24}
\LL_Y F(y) = c \frac{d^2}{dy^2} \int_0^\infty F(x) \log\biggl(
\frac{1}{\vert y \m x \vert} \biggr) \,dx + \frac{c F(0)}{y},
\end{equation}
\item[]\mbox{Caputo's form}\vspace*{-1pt}
\begin{equation}\label{A.25}
\LL_Y
F(y) = c \int_0^\infty F''(x) \log\biggl( \frac{1}{\vert y \m
x \vert}\biggr) \,dx,
\end{equation}
\end{itemize}
whenever $F \in C_b^2(\R_+)$ satisfies
%
%
\begin{equation}
\label{A.26}
F'(0+) = 0 \qquad\mbox{(normal reflection)}
\end{equation}
with $\vert F''(y) \vert= O(y^\gamma)$ as $y \rightarrow
\infty$ for some $\gamma< -1$ [as well as $\vert F(y) \vert=
O(y^\delta)$ as $y \rightarrow\infty$ for some $\delta< -1$ in
(\ref{A.24}) above].
\end{Proposition}
\begin{pf}
As in the proof of Proposition \ref{prop4} it is enough to
derive (\ref{A.23}). Using the same arguments as in (\ref{A.8}) and
(\ref{A.9}), we find that
%
%
\begin{eqnarray}
\label{A.27}
F(Y_t) &=& F(Y_0) - \int_0^t F'(Y_{s-}) \,dX_s\nonumber\\[-8pt]\\[-8pt]
&&{} + \sum_{0 < s \le t}
\bigl(F(Y_{s-} \p\Delta Y_s) - F(Y_{s-}) + F'(Y_{s-}) \Delta
X_s \bigr),\nonumber
\end{eqnarray}
where $(\int_0^t F'(Y_{s-}) \,dX_s)_{t \ge0}$ is a local
martingale. We can no longer claim that this process is a
martingale, however, we note from (\ref{A.22}) that $X_t = M_t +
A_t$ with
%
%
\begin{eqnarray}
\label{A.28}
\EE e^{i \lambda M_t} &=& \exp\biggl( t \int_{\vert x \vert\le1}
(e^{i \lambda x} - 1 - i \lambda x ) \nu_1(dx)\biggr), \vadjust{\goodbreak}\\[-3pt]
\label{A.29}
\EE e^{i \lambda A_t} &=& \exp\biggl( t
\int_{\vert x \vert> 1} ( e^{i \lambda x} - 1 )
\nu_1(dx) \biggr)
\end{eqnarray}
from where we see that the (L\'evy) process $M = (M_t)_{t \ge0}$ is
a martingale (whose L\'evy measure has bounded support) and the
bounded variation (L\'evy) process $A = (A_t)_{t \ge0}$ is given by
%
%
\begin{equation}
\label{A.30}
A_t = \sum_{0 < s \le t} \Delta X_s I(\vert\Delta X_s \vert> 1)
\end{equation}
for $t \ge0$. From (\ref{A.27})--(\ref{A.30}), we see that
%
%
\begin{eqnarray}
\label{A.31}
F(Y_t) &=& F(Y_0) - \int_0^t F'(Y_{s-}) \,dM_s \nonumber\\[-2pt]
&&{}
+
\sum_{0 < s \le t} \bigl(F(Y_{s-} \p\Delta Y_s) - F(Y_{s-})\nonumber\\[-2pt]
&&\hspace*{39.1pt}{}
+ F'(Y_{s-})
\Delta X_s I(\vert\Delta X_s \vert\le1) \bigr) \nonumber\\[-2pt]
&=& F(Y_0) - \int_0^t F'(Y_{s-}) \,dM_s \nonumber\\[-10pt]\\[-10pt]
&&{}
+ \sum_{0
< s \le t} \bigl( [ F(Y_{s-} - \Delta X_s) \m F(Y_{s-})\nonumber\\[-2pt]
&&\hspace*{42.8pt}{} \p F'(Y_{s-})
\Delta X_s I(\vert\Delta X_s \vert\le1) ] I(\Delta
X_s \le Y_{s-})\nonumber\\[-2pt]
&&\hspace*{39.4pt}{}
+ [F(0) \m F(Y_{s-})\nonumber\\[-2pt]
&&\hspace*{54.4pt}{} +
F'(Y_{s-}) \Delta X_s I(\vert\Delta X_s \vert\le1) ]
I(\Delta X_s > Y_{s-}) \bigr)\nonumber
\end{eqnarray}
using the same arguments as in (\ref{4.9}) above. Since $F'$ is
bounded and the L\'evy measure of $M$ has bounded support (implying
${\EE\sup_{ 0 \le s \le t} }\vert M_s \vert^q < \infty$ and hence
$\EE[M,M]^{q/2} < \infty$ for all $q>0$ by the BDG inequality) it
also follows by the BDG inequality (with $q=1$) that $(\int_0^t
F'(Y_{s-}) \,dX_s)_{t \ge0}$ is a~martingale. Taking $\EE_y$ on
both sides of (\ref{A.31}), where $\PP_{ y}$ denotes a probability
measure under which $Y_0=y$, and applying the compensation formula
(see, e.g., \cite{RY}, page 475) we find that
%
%
\begin{eqnarray}
\label{A.32}
&&\EE_y F(Y_t) \m F(y) \hspace*{-22pt}\nonumber\\[-2pt]
&&\qquad= \EE_y \biggl[\int_0^t ds \biggl(\int_{-
\infty}^{Y_s} [ F(Y_s \m x) \m F(Y_s) \p F'(Y_s) x I(
\vert x \vert\le1) ] \nu_1(dx) \hspace*{-22pt}\\[-2pt]
&&\qquad\quad\hspace*{63.7pt}{}
+
\int_{Y_s}^\infty[ F(0) \m F(Y_s) \p F'(Y_s) x
I(\vert x \vert\le1) ] \nu_1(dx) \biggr) \biggr]\hspace*{-22pt}\nonumber
\end{eqnarray}
for all $y>0$. The applicability of this formula (see, e.g.,
\cite{Ky}, page 97) follows from the facts that $\vert F(y) \vert\le
C (1 \p y^{\gamma+2})$, $\vert F'(y) \vert\le C (y \wedge
y^{\gamma+1})$ and $\vert F''(y) \vert\le C (1 \wedge y^\gamma)$
for all $y \ge0$ with some\vadjust{\goodbreak} $C>0$ so that the mean value theorem
yields the existence of $\xi_{s,x}^1 \in(Y_s,Y_s \p x)$,
$\xi_{s,x}^2 \in(Y_s \m x,Y_s))$ and $\eta_s \in(0,Y_s)$ such
that
%
%
\begin{eqnarray}
\label{A.33}\hspace*{13pt}
&&\EE_y \biggl[\int_0^t ds \biggl(\int_{-
\infty}^{Y_s} \bigl\vert F(Y_s \m x) \m F(Y_s) \p F'(Y_s) x I(
\vert x \vert\le1) \bigr\vert\nu_1(dx) \nonumber\\[1pt]
&&\qquad\quad\hspace*{30.1pt}{}
+
\int_{Y_s}^\infty\bigl\vert F(0) \m F(Y_s) \p F'(Y_s) x
I(\vert x \vert\le1) \bigr\vert\nu_1(dx) \biggr) \biggr] \nonumber\\[1pt]
&&\qquad\le\EE_y \biggl[\int_0^t ds \biggl( \int_1^\infty\vert
F(Y_s \p x) \m F(Y_s) \vert\nu_1(dx)\nonumber\\[1pt]
&&\qquad\quad\hspace*{48pt}{}
+
\int_0^1 \vert F(Y_s \p x) \m F(Y_s) \m F'(Y_s) x
\vert\nu_1(dx) \nonumber\\[1pt]
&&\qquad\quad\hspace*{48pt}{}
+ \int_0^{Y_s} \vert
F(Y_s \m x) \m F(Y_s) \p F'(Y_s) x \vert\nu_1(dx)
I(Y_s < 1) \nonumber\\[1pt]
&&\qquad\quad\hspace*{48pt}{}
+ \int_0^1 \vert F(Y_s \m x)
\m F(Y_s) \p F'(Y_s) x \vert\nu_1(dx) I(Y_s \ge1)
\nonumber\\[1pt]
&&\qquad\quad\hspace*{48pt}{}+ \int_1^{Y_s} \vert F(Y_s \m x) \m F(Y_s)
\vert\nu_1(dx) I(Y_s \ge1) \nonumber\\[1pt]
&&\qquad\quad\hspace*{48pt}{}+
\int_{Y_s}^1 \vert F(0) \m F(Y_s) \p F'(Y_s) x \vert
\nu_1(dx) I(Y_s < 1)\nonumber\\[1pt]
&&\qquad\quad\hspace*{164pt}{}+ \int_1^\infty
\vert F(0) \m F(Y_s) \vert\nu_1(dx) \biggr) \biggr] \nonumber\\[1pt]
&&\qquad\le\EE_y \biggl[\int_0^t ds \biggl( 2 C \int_1^\infty
\bigl(1 \p(Y_s \p x)^{\gamma+2} \bigr) \frac{c}{x^2} \,dx\nonumber\\[1pt]
&&\qquad\quad\hspace*{48pt}{} + \int_0^1
\frac{1}{2} \vert F''(\xi_{s,x}^1) \vert x^2 \frac{c}{x^2}
\,dx + 2 \int_0^1 \frac{1}{2} \vert F''(\xi_{s,x}^2)
\vert x^2 \frac{c}{x^2} \,dx\nonumber\\[1pt]
&&\qquad\quad\hspace*{48pt}{} + 2 C (1 \p Y_s^{\gamma+2})
\int_1^{Y_s} \frac{1}{x^2} \,dx\, I(Y_s \ge1) \nonumber\\[1pt]
&&\qquad\quad\hspace*{48pt}{}+
\int_{Y_s}^1 \bigl( \vert F'(\eta_s) \vert Y_s + \vert F'(Y_s)
\vert x \bigr) \frac{c}{x^2} \,dx\, I(Y_s < 1) \nonumber\\[1pt]
&&\qquad\quad\hspace*{48pt}\hspace*{77.5pt}{}+ \int_1^\infty\bigl(\vert F(0) \vert+ C ( 1 \p
Y_s^{\gamma+2} ) \bigr) \frac{c}{x^2} \,dx \biggr) \biggr] \nonumber\\[1pt]
&&\qquad\le c \EE_y \biggl[\int_0^t ds \biggl(4 C ( 1 \p Y_s^{
\gamma+2} ) \int_1^\infty\frac{1}{x^2} \,dx + 2 C \int_1^\infty
x^\gamma \,dx + \frac{3}{2} C \nonumber\\[1pt]
&&\qquad\quad\hspace*{52.1pt}{}
+ 2
C (1 \m Y_s) I(Y_s < 1) + \vert F(0) \vert+ C (1
\p Y_s^{\gamma+2}) \biggr) \biggr] \\
&&\qquad\le c \EE_y \biggl[
\int_0^t ds \biggl( \frac{17}{2} C + 5 C Y_s^{\gamma+2} -
\frac{2C}{\gamma+1} + \vert F(0) \vert\biggr) \biggr] \nonumber\\
&&\qquad\le c \biggl[\biggl(\frac{17}{2} C - \frac{2C}{\gamma+1} + \vert
F(0) \vert\biggr) t \nonumber\\
&&\qquad\quad\hspace*{9pt}{}+ \frac{5 C}{\gamma+3} t^{\gamma+3}
\EE_y(S_1 \m I_1)^{\gamma+2} \biggr] < \infty\nonumber
\end{eqnarray}
since $\gamma\p2 \in(0,1)$ and where we also use the scaling
property of $X$. (Note that without loss of generality we can assume
that $\gamma$ is close enough to $-1$ so that $\gamma\p2 > 0$.)
Dividing both sides of (\ref{A.32}) by $t$, letting $t \downarrow0$
and using the dominated convergence theorem, we get
%
%
\begin{eqnarray}
\label{A.34}
\LL_Y F(y) &=& \int_{-\infty}^y [ F(y \m x) \m F(y) \p F'(y)
x I(\vert x \vert\le1) ] \nu_1(dx) \nonumber\\
&&{}
+
[ F(0) \m F(y) ] \int_y^\infty\nu_1(dx) \\
&&{}+ F'(y) \int_y^\infty
x I(\vert x \vert\le1) \nu_1(dx)\nonumber
\end{eqnarray}
for all $y>0$. Splitting the integral over $(-\infty,y]$ into
integrals over $(-\infty,-y]$ and $[-y,y]$, noting that the third
term of the resulting integral over $(-\infty,-y]$ cancels with the
final term in (\ref{A.34}), it is easily seen using (\ref{A.2}) that
the expression on the right-hand side of (\ref{A.34}) coincides with
the expression on the right-hand side of (\ref{A.23}). This
completes the proof.
\end{pf}
\end{appendix}

%

%
\printaddresses


\begin{thebibliography}{27}

\bibitem{AS}
\begin{bbook}[mr]
\beditor{\bsnm{Abramowitz},~\bfnm{M.}\binits{M.}} \AND
  \beditor{\bsnm{Stegun},~\bfnm{I.~A.}\binits{I.~A.}}, eds.
(\byear{1992}).
\btitle{Handbook of Mathematical Functions with Formulas, Graphs, and
  Mathematical Tables}.
\bpublisher{Dover}, \baddress{New York}.
\bid{mr={1225604}}
\end{bbook}
\endbibitem

\bibitem{BDP}
\begin{barticle}[mr]
\bauthor{\bsnm{Bernyk},~\bfnm{Violetta}\binits{V.}},
  \bauthor{\bsnm{Dalang},~\bfnm{Robert~C.}\binits{R.~C.}} \AND
  \bauthor{\bsnm{Peskir},~\bfnm{Goran}\binits{G.}}
(\byear{2008}).
\btitle{The law of the supremum of a stable {L}\'evy process with no negative
  jumps}.
\bjournal{Ann. Probab.}
\bvolume{36}
\bpages{1777--1789}.
\bid{doi={10.1214/07-AOP376}, mr={2440923}}
\end{barticle}
\endbibitem

\bibitem{Be}
\begin{bbook}[mr]
\bauthor{\bsnm{Bertoin},~\bfnm{Jean}\binits{J.}}
(\byear{1996}).
\btitle{L\'evy Processes}.
\bseries{Cambridge Tracts in Mathematics}
\bvolume{121}.
\bpublisher{Cambridge Univ. Press}, \baddress{Cambridge}.
\bid{mr={1406564}}
\end{bbook}
\endbibitem

\bibitem{Bo}
\begin{barticle}[mr]
\bauthor{\bsnm{Boyce},~\bfnm{William~M.}\binits{W.~M.}}
(\byear{1970}).
\btitle{Stopping rules for selling bounds}.
\bjournal{Bell J. Econom. and Management Sci.}
\bvolume{1}
\bpages{27--53}.
\bid{mr={0285259}}
\end{barticle}
\endbibitem

\bibitem{Do}
\begin{barticle}[mr]
\bauthor{\bsnm{Doney},~\bfnm{R.~A.}\binits{R.~A.}}
(\byear{2008}).
\btitle{A note on the supremum of a stable process}.
\bjournal{Stochastics}
\bvolume{80}
\bpages{151--155}.
\bid{doi={10.1080/17442500701830399}, mr={2402160}}
\end{barticle}
\endbibitem

\bibitem{DP-1}
\begin{barticle}[mr]
\bauthor{\bparticle{du~}\bsnm{Toit},~\bfnm{J.}\binits{J.}} \AND
  \bauthor{\bsnm{Peskir},~\bfnm{G.}\binits{G.}}
(\byear{2007}).
\btitle{The trap of complacency in predicting the maximum}.
\bjournal{Ann. Probab.}
\bvolume{35}
\bpages{340--365}.
\bid{doi={10.1214/009117906000000638}, mr={2303953}}
\end{barticle}
\endbibitem

\bibitem{DP-2}
\begin{bincollection}[mr]
\bauthor{\bparticle{du~}\bsnm{Toit},~\bfnm{Jacques}\binits{J.}} \AND
  \bauthor{\bsnm{Peskir},~\bfnm{Goran}\binits{G.}}
(\byear{2008}).
\btitle{Predicting the time of the ultimate maximum for {B}rownian motion with
  drift}.
In \bbooktitle{Mathematical Control Theory and Finance}
\bpages{95--112}.
\bpublisher{Springer}, \baddress{Berlin}.
\bid{doi={10.1007/978-3-540-69532-5_6}, mr={2484106}}
\end{bincollection}
\endbibitem

\bibitem{DPS}
\begin{barticle}[mr]
\bauthor{\bparticle{du~}\bsnm{Toit},~\bfnm{J.}\binits{J.}},
  \bauthor{\bsnm{Peskir},~\bfnm{G.}\binits{G.}} \AND
  \bauthor{\bsnm{Shiryaev},~\bfnm{A.~N.}\binits{A.~N.}}
(\byear{2008}).
\btitle{Predicting the last zero of {B}rownian motion with drift}.
\bjournal{Stochastics}
\bvolume{80}
\bpages{229--245}.
\bid{doi={10.1080/17442500701840950}, mr={2402166}}
\end{barticle}
\endbibitem

\bibitem{GM}
\begin{barticle}[mr]
\bauthor{\bsnm{Gilbert},~\bfnm{John~P.}\binits{J.~P.}} \AND
  \bauthor{\bsnm{Mosteller},~\bfnm{Frederick}\binits{F.}}
(\byear{1966}).
\btitle{Recognizing the maximum of a sequence}.
\bjournal{J. Amer. Statist. Assoc.}
\bvolume{61}
\bpages{35--73}.
\bid{mr={0198637}}
\end{barticle}
\endbibitem

\bibitem{GR}
\begin{bbook}[vtex]
\bauthor{\bsnm{Gradshteyn},~\bfnm{I.~S.}\binits{I.~S.}} \AND
  \bauthor{\bsnm{Ryzhik},~\bfnm{I.~M.}\binits{I.~M.}}
(\byear{1994}).
\btitle{Table of Integrals, Series, and Products}.
\bpublisher{Academic Press}, \baddress{Boston, MA}.
\bid{mr={1243179}}
\end{bbook}
\endbibitem\vadjust{\goodbreak}

\bibitem{GPS}
\begin{barticle}[mr]
\bauthor{\bsnm{Graversen},~\bfnm{S.~E.}\binits{S.~E.}},
  \bauthor{\bsnm{Peskir},~\bfnm{G.}\binits{G.}} \AND
  \bauthor{\bsnm{Shiryaev},~\bfnm{A.~N.}\binits{A.~N.}}
(\byear{2001}).
\btitle{Stopping {B}rownian motion without anticipation as close as possible to
  its ultimate maximum}.
\bjournal{Theory Probab. Appl.}
\bvolume{45}
\bpages{41--50}.
\bid{doi={10.1137/S0040585X97978075}, mr={1810977}}
\bptnote{check year}
\end{barticle}
\endbibitem

\bibitem{GS}
\begin{barticle}[vtex]
\bauthor{\bsnm{Griffeath},~\bfnm{David}\binits{D.}} \AND
  \bauthor{\bsnm{Snell},~\bfnm{J.~Laurie}\binits{J.~L.}}
(\byear{1974}).
\btitle{Optimal stopping in the stock market}.
\bjournal{Ann. Probab.}
\bvolume{2}
\bpages{1--13}.
\bid{mr={0362766}}
\end{barticle}
\endbibitem

\bibitem{Ho}
\begin{bbook}[vtex]
\bauthor{\bsnm{Hochstadt},~\bfnm{Harry}\binits{H.}}
(\byear{1973}).
\btitle{Integral Equations}.
\bpublisher{Wiley}, \baddress{New York}.
\bid{mr={0390680}}
\end{bbook}
\endbibitem

\bibitem{Ka}
\begin{bincollection}[vtex]
\bauthor{\bsnm{Karlin},~\bfnm{Samuel}\binits{S.}}
(\byear{1962}).
\btitle{Stochastic models and optimal policy for selling an asset}.
In \bbooktitle{Studies in Applied Probability and Management Science}
\bpages{148--158}.
\bpublisher{Stanford Univ. Press}, \baddress{Stanford, CA}.
\bid{mr={0137594}}
\end{bincollection}
\endbibitem

\bibitem{Ky}
\begin{bbook}[vtex]
\bauthor{\bsnm{Kyprianou},~\bfnm{Andreas~E.}\binits{A.~E.}}
(\byear{2006}).
\btitle{Introductory Lectures on Fluctuations of {L}\'evy Processes with
  Applications}.
\bpublisher{Springer}, \baddress{Berlin}.
\bid{mr={2250061}}
\end{bbook}
\endbibitem

\bibitem{Pa}
\begin{barticle}[mr]
\bauthor{\bsnm{Patie},~\bfnm{P.}\binits{P.}}
(\byear{2009}).
\btitle{A few remarks on the supremum of stable processes}.
\bjournal{Statist. Probab. Lett.}
\bvolume{79}
\bpages{1125--1128}.
\bid{doi={10.1016/j.spl.2009.01.001}, mr={2510779}}
\end{barticle}
\endbibitem

\bibitem{Pe}
\begin{barticle}[mr]
\bauthor{\bsnm{Pedersen},~\bfnm{Jesper~Lund}\binits{J.~L.}}
(\byear{2003}).
\btitle{Optimal prediction of the ultimate maximum of {B}rownian motion}.
\bjournal{Stoch. Stoch. Rep.}
\bvolume{75}
\bpages{205--219}.
\bid{doi={10.1080/1045112031000118994}, mr={1994906}}
\end{barticle}
\endbibitem

\bibitem{PS}
\begin{bbook}[vtex]
\bauthor{\bsnm{Peskir},~\bfnm{Goran}\binits{G.}} \AND
  \bauthor{\bsnm{Shiryaev},~\bfnm{Albert}\binits{A.}}
(\byear{2006}).
\btitle{Optimal Stopping and Free-Boundary Problems}.
\bpublisher{Birkh\"auser}, \baddress{Basel}.
\bid{mr={2256030}}
\end{bbook}
\endbibitem

\bibitem{Pr}
\begin{bbook}[vtex]
\bauthor{\bsnm{Protter},~\bfnm{Philip~E.}\binits{P.~E.}}
(\byear{2004}).
\btitle{Stochastic Integration and Differential Equations}.
\bpublisher{Springer}, \baddress{Berlin}.
\end{bbook}
\endbibitem

\bibitem{RY}
\begin{bbook}[mr]
\bauthor{\bsnm{Revuz},~\bfnm{Daniel}\binits{D.}} \AND
  \bauthor{\bsnm{Yor},~\bfnm{Marc}\binits{M.}}
(\byear{1999}).
\btitle{Continuous Martingales and {B}rownian Motion},
\bedition{3rd} ed.
\bseries{Grundlehren der Mathematischen Wissenschaften [Fundamental Principles
  of Mathematical Sciences]}
\bvolume{293}.
\bpublisher{Springer}, \baddress{Berlin}.
\bid{mr={1725357}}
\end{bbook}
\endbibitem

\bibitem{Sa}
\begin{bbook}[vtex]
\bauthor{\bsnm{Sato},~\bfnm{Ken-iti}\binits{K.-i.}}
(\byear{1999}).
\btitle{L\'evy Processes and Infinitely Divisible Distributions}.
\bseries{Cambridge Studies in Advanced Mathematics}
\bvolume{68}.
\bpublisher{Cambridge Univ. Press}, \baddress{Cambridge}.
\bid{mr={1739520}}
\end{bbook}
\endbibitem

\bibitem{Sh-1}
\begin{bincollection}[vtex]
\bauthor{\bsnm{Shiryaev},~\bfnm{Albert~N.}\binits{A.~N.}}
(\byear{2002}).
\btitle{Quickest detection problems in the technical analysis of the financial
  data}.
In \bbooktitle{Mathematical Finance---{B}achelier {C}ongress, 2000 ({P}aris)}
\bpages{487--521}.
\bpublisher{Springer}, \baddress{Berlin}.
\bid{mr={1960576}}
\end{bincollection}
\endbibitem

\bibitem{Sh-2}
\begin{barticle}[vtex]
\bauthor{\bsnm{Shiryaev},~\bfnm{A.~N.}\binits{A.~N.}}
(\byear{2009}).
\btitle{On conditional-extremal problems of the quickest detection of
  nonpredictable times of the observable Brownian motion}.
\bjournal{Theory Probab. Appl.}
\bvolume{53}
\bpages{663--678}.%
\end{barticle}%
\endbibitem%

\bibitem{Ur}
\begin{barticle}[mr]
\bauthor{\bsnm{Urusov},~\bfnm{M.~A.}\binits{M.~A.}}
(\byear{2005}).
\btitle{On a property of the time of attaining the maximum by {B}rownian motion
  and some optimal stopping problems}.
\bjournal{Theory Probab. Appl.}
\bvolume{49}
\bpages{169--176}.
\bid{doi={10.1137/S0040585X97980956}, mr={2141339}}
\bptnote{check year}
\end{barticle}
\endbibitem

\bibitem{Wa}
\begin{barticle}[mr]
\bauthor{\bsnm{Watanabe},~\bfnm{Shinzo}\binits{S.}}
(\byear{1962}).
\btitle{On stable processes with boundary conditions}.
\bjournal{J. Math. Soc. Japan}
\bvolume{14}
\bpages{170--198}.
\bid{mr={0144387}}
\end{barticle}
\endbibitem

\bibitem{Wi}
\begin{bbook}[vtex]
\bauthor{\bsnm{Williams},~\bfnm{David}\binits{D.}}
(\byear{1991}).
\btitle{Probability with Martingales}.
\bpublisher{Cambridge Univ. Press}, \baddress{Cambridge}.
\bid{mr={1155402}}
\end{bbook}
\endbibitem

\end{thebibliography}
\end{document}